\documentclass[11pt]{amsart}

\usepackage{amsfonts, amstext, amsmath, amsthm, amscd, amssymb}
\usepackage{epsfig, graphics, psfrag}
\usepackage{color}

\renewcommand{\setminus}{{\smallsetminus}}

\newcommand{\cross}{{\times}} 
\newcommand{\bdy}{{\partial}} 

\newcommand{\half}{{\frac{1}{2}}}
\newcommand{\abs}[1]{{\left\vert #1 \right\vert}}
\renewcommand{\inf}[1]{{\mathrm{inf} \left\{ #1 \right\} }}

\newcommand{\no}{{\noindent}}
\newcommand{\e}{{\epsilon}}
\newcommand{\kte}{{k_{t,\e}}}
\newcommand{\fte}{{f_{t,\e}}}
\newcommand{\fpte}{{f'_{t,\e}}}
\newcommand{\fppte}{{f''_{t,\e}}}
\newcommand{\gte}{{g_{t,\e}}}
\newcommand{\gpte}{{g'_{t,\e}}}
\newcommand{\tlim}{{t_{\rm lim}}}
\newcommand{\elim}{{\e_{\rm lim}}}

\newcommand{\GA}{{\mathbb{G}_A}}
\newcommand{\GB}{{\mathbb{G}_B}}
\newcommand{\GRA}{{\mathbb{G}'_A}}
\newcommand{\GRB}{{\mathbb{G}'_B}}
\newcommand{\vbigon}{{v_{\rm bigon}}}
\newcommand{\vlarge}{{v_{\rm n\minus gon}}}
\newcommand{\eshort}{{e_{\rm short}}}
\newcommand{\elong}{{e_{\rm long}}}
\newcommand{\ershort}{{e'_{\rm short}}}
\newcommand{\erlong}{{e'_{\rm long}}}
\newcommand{\lmin}{{\ell_{\rm min}}}
\newcommand{\twist}{{\mathrm{tw }(D)}}
\newcommand{\twistcountry}{{\mathrm{tw }(N)}}

\DeclareMathSymbol{\minus} {\mathord}{operators}{"2D} 
\newcommand{\vol}{{\rm vol}}
\newcommand{\area}{{\rm area}}

\renewcommand{\phi}{{\varphi}}

\newcommand{\RR}{{\mathbb{R}}}

\theoremstyle{plain}
\newtheorem{theorem}{Theorem}[section]
\newtheorem{corollary}[theorem]{Corollary}
\newtheorem{lemma}[theorem]{Lemma}
\newtheorem{prop}[theorem]{Proposition}

\newtheorem*{no-num-theorem}{Theorem}

\theoremstyle{definition}
\newtheorem{define}[theorem]{Definition}
\newtheorem*{remark}{Remark}

\begin{document}

\title[Dehn filling, Volume, and the Jones polynomial]{Dehn filling, Volume, and \\ the Jones polynomial}
\author[D.\ Futer]{David Futer}
\author[E.\ Kalfagianni]{Efstratia Kalfagianni}
\author[J.\ Purcell]{Jessica S. Purcell}
\date{12/18/2006}

\thanks{The first author was supported in part by  NSF--RTG grant
	DMS--0353717. The second author was supported in part by NSF grant
	DMS--0306995  and NSF--FRG grant DMS-0456155.  The third author was
	supported in part by NSF grant DMS-0704359.} 


\address{David Futer, Mathematics Department, Michigan State
  University, East Lansing, MI 48824}
\email{dfuter@math.msu.edu}

\address{Efstratia Kalfagianni, Mathematics Department, Michigan State
  University, East Lansing, MI 48824}
\email{kalfagia@math.msu.edu}

\address{Jessica S. Purcell, Department of Mathematics, Brigham Young
  University, Provo, UT 84602} \email{jpurcell@math.byu.edu}

\begin{abstract}
  Given a hyperbolic 3--manifold with torus boundary, we bound the
  change in volume under a Dehn filling where all slopes have length
  at least $2\pi$.  This result is applied to give explicit diagrammatic 
  bounds on the volumes of many knots and links, as well as
  their Dehn fillings and branched covers.  Finally, we use this 
  result to bound the volumes
  of knots in terms of the coefficients of their Jones polynomials.
\end{abstract}

\maketitle

\section{Introduction \label{sec:intro}}
It is well--known that the volumes of hyperbolic 3--manifolds form a closed, 
well--ordered subset of $\RR$ \cite{thurston-notes}.  However, 3--manifolds 
are often described combinatorially, and it remains hard to translate the 
combinatorial data into explicit information on volume. 
 In this paper, we prove results that bound the volumes of a large class
of manifolds with purely combinatorial descriptions.  

There are other recent theorems relating volumes to combinatorial
data.  Brock and Souto have proved that the volume of a hyperbolic
3--manifold is coarsely determined by the complexity of a Heegaard
splitting \cite{brock-souto}.  Costan\-tino and Thurston have related
volume to the complexity of a shadow complex \cite{costantino-thurston}.
Despite the general power of these theorems, the constants that bound
volume from below remain mysterious.

This paper provides explicit and readily applicable estimates on the
volume of hyperbolic manifolds obtained using Dehn filling.  We apply
these estimates to a large class of knot and link complements,
obtaining bounds on their volume based purely on the combinatorics of
a diagram of the link.  We then use these results to relate the volume of a
large class of knots to the coefficients of the Jones polynomial.

The volume conjecture \cite{kashaev:volume-conj, murakami-squared}
asserts that the volume of hyperbolic knots is determined by certain
asymptotics of the Jones polynomial and its relatives.  At the same
time, a wealth of experimental evidence suggests a direct correlation
between the coefficients of the Jones polynomial and the volume of
hyperbolic knots. For example, Champanerkar, Kofman, and Patterson
have computed the Jones polynomials of all the hyperbolic knots whose
complements can be decomposed into seven or fewer ideal tetrahedra
\cite{ckp}.  Although some of these Jones polynomials have large
spans, their non-zero coefficients have small values, suggesting a
relationship between small volume and small coefficients.  Dasbach and
Lin have proved that such a connection does in fact exist for
alternating links \cite{dasbach-lin:volumeish}; our results extend
this relationship to many non-alternating links.

\medskip

\subsection{Volume change}

Given a 3--manifold $M$ with $k$ torus boundary components, we use the
following standard terminology.  For the $i$-th torus $T_i$, let $s_i$
be a \emph{slope} on $T_i$, that is, an isotopy class of simple closed
curves.  Let $M(s_1, \dots, s_k)$ denote the manifold obtained by Dehn
filling $M$  along the slopes $s_1$, \dots, $s_k$.

When $M$ is hyperbolic, each torus boundary component of $M$
corresponds to a cusp.  Taking a maximal disjoint horoball
neighborhood about the cusps, each torus $T_i$ inherits a Euclidean
structure, well--defined up to similarity.  The slope $s_i$ can then
be given a geodesic representative.  We define the \emph{slope length}
of $s_i$ to be the length of this geodesic representative.  Note that when
$k>1$, this definition of slope length depends on the choice of
maximal horoball neighborhood.  

\begin{theorem}\label{thm:vol-change}
  Let $M$ be a complete, finite--volume hyperbolic manifold with
  cusps.  Suppose $C_1, \dots, C_k$ are disjoint horoball
  neighborhoods of some subset of the cusps.  Let $s_1, \dots, s_k$ be
  slopes on $\partial C_1, \dots, \partial C_k$, each with length
  greater than $2\pi$.  Denote the minimal slope length by $\lmin$.
  If $M(s_1, \dots, s_k)$ satisfies the geometrization conjecture,
  then it is a hyperbolic manifold, and
  $$ \vol(M(s_1, \dots, s_k)) \: \geq \:
  \left(1-\left(\frac{2\pi}{\lmin}\right)^2\right)^{3/2} \vol(M).$$
\end{theorem}

Note that when at least one cusp of $M$ is left
unfilled, the manifold $M(s_1, \dots, s_k)$ is Haken, and thus
satisfies geometrization by Thurston's theorem \cite{thurston-survey}.
In the general case, the hyperbolicity of $M(s_1, \dots, s_k)$ would
follow from Perelman's work \cite{perelman02, perelman03}.

Theorem \ref{thm:vol-change} should be compared with other known
results.  Neumann and Zagier have found asymptotic changes in volume
under Dehn filling as slope lengths become long \cite{neumann-zagier}.
They show that the change in volume is asymptotically of order
$O(1/\lmin^2)$.  Although Theorem \ref{thm:vol-change} was not meant
to analyze the asymptotic behavior of volume, it also gives an
$O(1/\lmin^2)$ estimate. However, our constants are not sharp. See
Section \ref{sec:numerics} for a more detailed discussion of the
sharpness and asymptotic behavior of our estimate.

Hodgson and Kerckhoff have also found bounds on volume change under
Dehn filling, provided that the filling is obtained via cone
deformation \cite{hk:univ}.  They show that if the normalized slope
length is at least 7.515, the cone deformation exists and their volume
estimates apply. However, the normalized slope length is typically
much smaller than the actual slope length.  Thus Theorem
\ref{thm:vol-change} applies in many more cases than their results.


\subsection{Twist number and volumes} \label{sec:volume-results}

We will apply Theorem \ref{thm:vol-change} to link complements in
$S^3$.  Consider a diagram of a knot or link $K$ as a 4--valent graph
in the plane, with over--under crossing information associated to each
vertex.  A bigon region is a region of the graph bounded by only two
edges.  A \emph{twist region} of a diagram consists of maximal
collections of bigon regions arranged end to end.  A single crossing
adjacent to no bigons is also a twist region.  Let $D(K)$ denote the
diagram of $K$.  We denote the number of twist regions in a diagram by
$\twist$.

\begin{figure}[h]
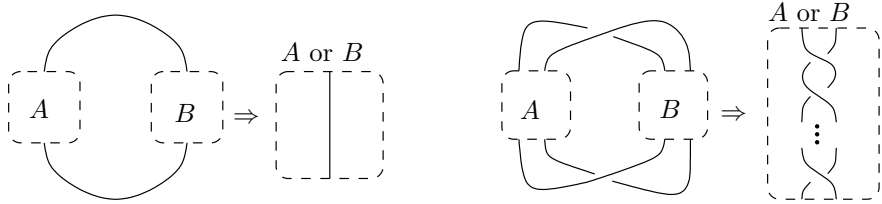

\begin{center}
\input{Fig_prime.pstex_t}
\hspace{.5in}
\input{Fig_twred.pstex_t}
\end{center}
\caption{Left:  A prime diagram.  Right:  A twist reduced diagram.}
\label{fig:prime-twred}
\end{figure}

Our statements concern the number of twist regions of a diagram.  We
rule out extraneous twist regions by requiring our diagram to be
reduced in the sense of the following two definitions, illustrated in
Figure \ref{fig:prime-twred}.

First, we require the diagram to be \emph{prime}.  That is, any simple
closed curve which meets two edges of the diagram transversely must
bound a region of the diagram with no crossings.

Second, we require the diagram to be \emph{twist--reduced}.  That is,
if any simple closed curve meets the diagram transversely in four
edges, with two points of intersection adjacent to one crossing and
the other two adjacent to another crossing, then that simple closed
curve must bound a (possibly empty) collection of bigons arranged end
to end between the crossings.

In the remainder of this paper, we will implicitly assume that all
link diagrams are connected, and that a diagram is alternating within
each twist region.

\begin{theorem}\label{thm:link-volume}
  Let $K \subset S^3$ be a link with a prime, twist--reduced diagram
  $D(K)$. Assume that $D(K)$ has $\twist \geq 2$ twist regions, and
  that each region contains at least $7$ crossings. Then $K$ is a
  hyperbolic link satisfying
  $$0.70735 \; (\twist - 1) \; < \; \vol(S^3 \setminus K) \; < \; 10\,
  v_3 \, (\twist - 1), $$
  where $v_3 \approx 1.0149$ is the volume of
  a regular ideal tetrahedron.
\end{theorem}

The upper bound on volume is due to Agol and D.\ Thurston
\cite[Appendix]{lackenby:alt-volume}, improving an earlier estimate by
Lackenby \cite{lackenby:alt-volume}. For alternating diagrams, Agol,
Storm, and W.\ Thurston \cite{ast-guts} have proved a sharper lower
bound of $1.83(\twist - 2)$, again improving an earlier estimate by
Lackenby \cite{lackenby:alt-volume}.  Theorem \ref{thm:link-volume} is
also an improvement of a recent theorem of Purcell
\cite{purcell:volume}.  A linear lower bound was also obtained in that
paper, but the results applied only to links with significantly more
crossings per twist region.

Theorem \ref{thm:vol-change} also leads to lower bounds on the volumes
of Dehn fillings of link complements in $S^3$ and branched coverings
of $S^3$ over links.  For example, combining Theorem
\ref{thm:vol-change} with the orbifold theorem \cite{blp, chk} and a
result of Adams on the waist size of knots
\cite{adams:waist}
yields the following result. 

\begin{theorem}\label{thm:branched} For a hyperbolic knot $K$ in $S^3$ and an
  integer $p>0$, let $M_p$ denote the $p$--fold cyclic cover of $S^3$
  branched over $K$. If  $p \geq 7$,
 then $M_p$ is hyperbolic, and
 $$
 \left(1-\frac{4\pi^2}{p^2}\right)^{3/2}
 \vol(S^3\setminus K) \: \leq \: {\vol(M_p)\over p}\: < \:
 \vol(S^3\setminus K).$$
\end{theorem}

For further applications and discussion, including a sharper version of Theorem \ref{thm:branched}, 
we refer the reader to Section 3.


\subsection{Twist number and Jones polynomials}

\begin{figure}[ht]
  \input{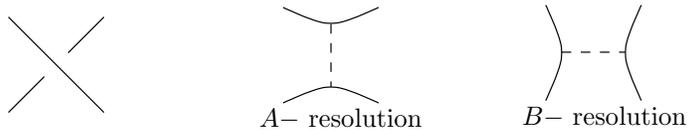} 
\caption{A crossing and its $A$--, $B$--resolutions. The dashed lines
  show the edges of the graphs $\GA$, $\GB$ corresponding to the
  crossing.} 
\label{fig:crossing-resolutions}
\end{figure}

Let $D$ be a link diagram, and $x$ a crossing of $D$. Associated to
$D$ and $x$ are two link diagrams, each with one fewer crossing than
$D$, called the \emph{$A$--resolution} and \emph{$B$--resolution} of
the crossing.  See Figure \ref{fig:crossing-resolutions}.  Starting
with any $D$, let $s_A(D)$ (resp. $s_B(D)$) denote the crossing--free
diagram obtained by applying the $A$--resolution (resp.
$B$--resolution) to all the crossings of $D$.
\smallskip

\vspace{-1ex}


\begin{define} \label{def:ab-graphs} Given a link diagram $D$ we
  obtain graphs $\GA$, $\GB$ as follows. The vertices of $\GA$ are in
  one-to-one correspondence with the components of $s_A(D)$.  Every
  crossing of $D$ gives rise to two arcs of the $A$--resolution.
  These will each be associated with a component of $s_A(D)$, and thus
  correspond to a vertex of $\GA$.  Add an edge to $\GA$ connecting
  these two vertices for each crossing of $D$, as in Figure
  \ref{fig:crossing-resolutions}.  We will refer to $\GA$ as the
  $A$--graph associated to $D$.  In a similar manner, construct the
  $B$--graph $\GB$ by considering components of $s_B(D)$.
  
  A link diagram $D$ is called $A$--adequate (resp. $B$--adequate) if
  the graph $\GA$ (resp. $\GB$) contains no loops (i.e. edges with
  both of their endpoints on the same vertex).  The diagram $D$ is
  called adequate if it is both $A$--adequate and $B$--adequate. A
  link is called adequate if it admits an adequate diagram.
\end{define}

The class of adequate links includes all alternating links and all
$n$--string parallels of alternating links, as well as most pretzel
knots and links and most arborescent links. For more information, see
for example the paper of Lickorish and Thistlethwaite
\cite{lick-thistle}.

For any link $K \subset S^3$, let
$$J_K(t)= \alpha t^n+ \beta t^{n-1}+ \ldots + \beta' t^{s+1}+ \alpha'
t^s$$ 
\no denote the Jones polynomial of $K$, so that $n$ (resp. $s$)
is the highest (resp. lowest) power in $t$. We will always denote the
second and next-to-last coefficients of $J_K(t)$ by $\beta$ and
$\beta'$, respectively.

\begin{theorem}\label{thm:twist-jones}
  Let $K$ be a link in $S^3$ with an adequate diagram $D(K)$, such
  that every twist region of $D(K)$ contains at least $3$ crossings.
  Then
  $$\frac{1}{3} \; \twist +1 \; \leq \; \abs{\beta} + \abs{\beta'} \; \leq \; 2
  \, \twist.$$
\end{theorem}

By putting together Theorem \ref{thm:link-volume} and Theorem
\ref{thm:twist-jones}, we obtain the following result relating the
volume and the Jones polynomial of a hyperbolic link.

\begin{corollary}\label{cor:jones-volume}
  Let $K \subset S^3$ be a link with a prime, twist--reduced, adequate
  diagram $D(K)$. Assume that $D(K)$ has $\twist \geq 2$ twist
  regions, and that each region contains at least $7$ crossings. Then
  $K$ is a hyperbolic link, satisfying
  $$0.35367 \, (\abs{\beta} + \abs{\beta'} - 2) \; < \; \vol(S^3
  \setminus K) \; < \; 30\, v_3 \, (\abs{\beta} + \abs{\beta'} - 1).
  $$
  Here, $\beta$ and $\beta'$ are the second and next-to-last
  coefficients of the Jones polynomial of $K$, and $v_3 \approx
  1.0149$ is the volume of a regular ideal tetrahedron.
\end{corollary}

Dasbach and Lin \cite{dasbach-lin:volumeish} showed that the twist
number of a twist--reduced alternating diagram is exactly $\abs{\beta}
+ \abs{\beta'}$.  Combined with work of Lackenby
\cite{lackenby:alt-volume}, this led to two--sided bounds on the
volume of alternating links in terms of these coefficients of the
Jones polynomial.  Theorem \ref{thm:twist-jones} and Corollary
\ref{cor:jones-volume} extend these results into the realm of
non-alternating links.  \smallskip

\subsection{Organization of the paper}

In Section \ref{sec:vol-change}, we prove Theorem \ref{thm:vol-change}
and provide some experimental data. The proof of \ref{thm:vol-change}
requires a careful analysis of the properties of solutions to certain
differential equations; due to their technical nature, these details
are postponed until Section \ref{sec:diffeq}.  In Section
\ref{sec:knots}, we apply Theorem \ref{thm:vol-change} to knots and
links, their Dehn fillings, and their brached covers.  In particular,
we prove Theorem \ref{thm:link-volume} and several other applications.
In Section \ref{sec:jones}, we relate the twist number of a diagram to
the Jones polynomial, proving Theorem \ref{thm:twist-jones}.

\subsection{Acknowledgements}

We thank Marc Lackenby for pointing us in the right direction with
differential equation arguments in the proof of Theorem
\ref{thm:neg-curved-solid-torus}. We thank Nathan Dunfield for helping
us set up the numerical experiments to check the sharpness of our
volume estimate. Finally, we are grateful to Lawrence Roberts, Peter
Storm, and Xiao\-dong Wang for their helpful suggestions.

\section{Volume change under filling \label{sec:vol-change}}

In this section, we prove Theorem \ref{thm:vol-change}, by employing
the following strategy.  For every cusp of $M$ that we need to fill,
we will explicitly construct a negatively curved solid torus,
following the proof of Gromov and Thurston's $2\pi$--theorem
\cite{bleiler-hodgson}. When we sew in these solid tori, we obtain a
negatively curved Riemannian metric on $M(s_1, \ldots, s_k)$. Then, we
will use a theorem of Boland, Connell, and Souto
\cite{boland-connell-souto} to compare the volume of this metric with
the true hyperbolic volume of the filled manifold.

This strategy is similar to that of Agol in \cite{agol:drilling}.
However, while Agol starts with closed hyperbolic manifolds and
constructs negatively curved metrics on cusped ones, we begin with
cusped hyperbolic manifolds and construct negatively curved metrics on
their Dehn fillings.

\subsection{Negatively curved metrics on a solid torus}\label{sec:solid-torus}
Our main tool in the proof of Theorem \ref{thm:vol-change} is the
following result, inspired by Cooper and Lackenby \cite[Proposition
3.1]{cooper-lackenby}. To simplify exposition, we define a function
$$h(x) := 1- \left( \frac{2\pi}{x} \right)^2.$$



\begin{theorem}\label{thm:neg-curved-solid-torus}
  Let $V$ be a solid torus.  Assume that $\bdy V$ carries a Euclidean metric,
in which the Euclidean geodesic representing a meridian has length
$\ell_1 > 2\pi$. Then, for any constant
$\zeta \in (0,1)$, there exists a
  smooth Riemannian metric $\tau$ on $V$, with the following
  properties:
\begin{itemize}
\item[(a)] On a collar neighborhood of $\bdy V$, $\tau$ is a
  hyperbolic metric, whose restriction to $\bdy V$ is the prescribed flat metric.
\item[(b)] The sectional curvatures of $\tau$ are bounded above by
$- \zeta\, h(\ell_1)$.
\item[(c)] The volume of $V$ in this metric is at least $\half \zeta\,
  \area(\bdy V)$.
\end{itemize}
\end{theorem}

\begin{proof}
  
  

Following Bleiler and Hodgson's proof of the $2\pi$ Theorem \cite{bleiler-hodgson}, we will explicitly construct a metric on $\widetilde{V}$, the universal cover of $V$. First, give $\widetilde{V}$ cylindrical coordinates
$(r, \mu, \lambda)$, where $r\leq 0$ is the radial distance measured
\emph{outward} from $\bdy \widetilde{V}$, $0\leq \mu\leq 1$ is measured around each meridional circle, 
and $-\infty < \lambda < \infty$ is measured in the londitudinal
direction, perpendicular to $\mu$. We normalize the coordinates so that the generator of the deck transformation group on $\widetilde{V}$ changes the $\lambda$ coordinate by $1$. 

The Riemannian metric on $\widetilde{V}$ is given by \begin{equation}
ds^2 = dr^2 + \left(f(r)\right)^2 d\mu^2 + \left(g(r)\right)^2
d\lambda^2,
\label{eqn:riemann-metric}
\end{equation}
where $f$ and $g$ are smooth functions that 
we will construct in the course of the proof. In order to obtain the prescribed Euclidean metric on $\bdy \widetilde{V}$, we must set $f(0) = \ell_1$ and $g(0) = \ell_2$, where $\ell_2 := \area(\bdy V) / \ell_1$. 

With this metric, the deck transformation group on $\widetilde{V}$ is generated by the isometry
$$(r, \mu, \lambda) \mapsto (r, \mu+\theta, \lambda+1),$$
where the
\emph{shearing factor} $\theta \in [0,1)$ is chosen so that the
fundamental domain of $\bdy V$ becomes a parallelogram of the correct
shape. See Figure \ref{fig:tuna-can}.  The metric on $\widetilde{V}$
descends to give a smooth metric on $V$, and the coordinates $(r, \mu,
\lambda)$ give local cylindrical coordinates on $V$.

\begin{figure}[h]
  \input{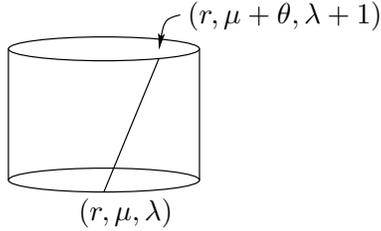}
  \caption{The fundamental domain for the action of the deck transformation group on $\widetilde {V}$.}
\label{fig:tuna-can}
\end{figure}

In order to give conclusions (a)--(c) of the theorem, the functions $f$ and $g$ must
satisfy several conditions:
\begin{itemize}
\item $f$ and $g$ must give a hyperbolic metric near $\bdy V$, such
  that the induced metric on $\bdy V$
  gives a Euclidean torus with the right shape.
  In other words, we must have $f(r) = \ell_1 e^r$ and $g(r)
  = \ell_2 e^r$ near $r=0$.
  
\item 
  In order to be nonsingular, the metric must have a
  cone angle of $2\pi$ along the core, i.e., at the points $r=r_0$
  such that $f(r_0)=0$.
Bleiler and Hodgson computed that this cone angle is exactly $f'(r_0)$.
  Thus we need to ensure $f'(r_0)=2\pi$.
  
\item Bleiler and Hodgson computed that the sectional curvatures are
  all convex combinations of:
  $$\kappa_{12} = -\frac{f''}{f}\, , \quad
  \kappa_{13} = -\frac{g''}{g}\, , \quad
  \kappa_{23} = - \frac{f' \cdot g'}{f \cdot g}\, .$$
  To ensure they are all bounded above by $-\zeta h(\ell_1)$, we
  ensure that each of these quantities is bounded.

\item The volume of $V$ is given by $\int_{r_0}^0 fg\,dr$.  For
  the volume estimate, we ensure this quantity is bounded below.
\end{itemize}


With these requirements in mind, we can begin to construct $f$ and
$g$. 
Basically, we construct both functions so that the curvature
estimate will be automatically true, and show that the other
conditions follow.  Roughly, we would like to fix a value $t>0$ and
define $f$ by a differential equation $f''/f = t$, and $g$ by $f'g'/fg
= t$.  This would imply that all curvatures are bounded above by $-t$.
However, this simple definition will not give a smooth hyperbolic
metric near $\bdy V$.  Thus we introduce smooth bump functions.

\begin{figure}[h]
\psfrag{e1}{$-\e$}
\psfrag{e2}{$-\e/2$}
\psfrag{r}{$r$}
\psfrag{t}{$t$}
\psfrag{1}{$1$}
\psfrag{0}{$0$}
\begin{center}
\includegraphics{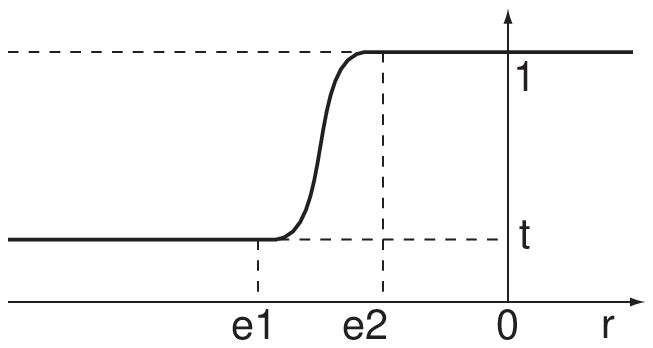}
\caption{The bump function $\kte(r)$.}
\label{fig:bump-function}
\end{center}
\end{figure}

For $\epsilon>0$ and $0 < t < 1$, let $\kte(r)$ be the smooth
bump function defined as follows:
$\kte(r) = t$ if $r\leq -\epsilon$, $\kte(r) = 1$ if $r\geq
-\epsilon/2$.  For $r$ between $-\epsilon$ and $-\epsilon/2$,
$\kte(r)$ is smooth and strictly increasing.
See Figure \ref{fig:bump-function} for a typical graph. We also extend
the definition of $\kte$ to $\epsilon=0$, obtaining a step function:
$$k_{t,0}(r) \; := \;  \lim_{\epsilon \to 0^+} \, \kte(r) \; = \; 
\left\{ 
\begin{array}{cl}
t &  \mbox{ if} \quad r < 0, \\
1 & \mbox{ if} \quad r \geq 0.
\end{array}
\right.
$$
Note that $k$ is continuous in the three variables $(r, t, \epsilon)$
for $\epsilon>0$.

For any $\e \geq 0$ and $t \in (0,1)$, define a function $\fte$
according to the differential equation
\begin{equation}
\fppte(r) = \kte(r) \fte(r),
\label{eqn:def-f}
\end{equation}
with initial conditions $\fte(0)= \ell_1 = \fpte(0)$.  
When $\epsilon > 0$ (and $k$ is continuous), the existence and
uniqueness of the solution $\fte$ is a standard result in differential
equations (see for example \cite{hsieh-sibuya}). When $\epsilon =0$,
the equation can be solved explicitly; $f_{t,0}$ is a $C^1$ function
that satisfies (\ref{eqn:def-f}) for all $r \neq 0$. (See equation
(\ref{eq:f0-explicit}) for the exact formula.)

In Section \ref{sec:diffeq}, we prove that the family of functions
$\fte$ has a number of nice properties. In particular, by Theorem
\ref{thm:f-regular}, $\fte(r)$ depends continuously and uniformly on
the parameters $t$ and $\e$. (When $\e>0$ and $\kte(r)$ is continuous,
this is a standard result in differential equations; when $\e \to 0$
and $\kte$ becomes discontinuous, this takes some work.)

Given $\fte$, we define 
$\gte$ according to the differential equation
\begin{equation}
\frac{\gpte(r)}{\gte(r)} \: := \:
\kte(r) \, \frac{\fte(r)}{\fpte(r)} \, ,
\label{eqn:def-g}
\end{equation}
with initial condition $\gte(0) = \ell_2$.  
Note that by Lemma \ref{lemma:f-monotonic}(c), $\fpte(r)>0$ for all
$r$, so the right-hand side is always well-defined. 
Because the left-hand side of (\ref{eqn:def-g}) is merely the
derivative of $\ln \gte(r)$,
the existence and uniqueness of solutions follows immediately by
integration.

\smallskip

Before we delve deeper into the properties of $\fte$ and $\gte$, a
roadmap is in order. By Lemma \ref{lemma:root-existence}, we know that
$\fte$ has a unique root $r_0 <0$. On the interval $[r_0, -\e]$,
$\fte$ and $\gte$ will have the form
$$\fte(r) = a \sinh(\sqrt{t}(r-r_0)), \quad \quad \gte(r) = b
\cosh(\sqrt{t}(r-r_0)), $$
for constants $a,b$ that depend on $t$ and $\e$. Qualitatively, this
means that the metric defined by equation $(\ref{eqn:riemann-metric})$
realizes the inner part of the solid torus as a rescaled hyperbolic
tube, with constant curvature $-t$ and a cone angle of $a\sqrt{t}$
along the core.

We will show that when $t = h(\ell_1)$ and $\e=0$, the cone angle is
exactly $2\pi$, and we get a non-singular tube of constant curvature
$-h(\ell_1)$. Furthermore, the volume of this metric is exactly $\half
\ell_1 \ell_2$. These values are certainly enough to satisfy
conditions (b) and (c) of the theorem. However, because $k_{t,0}$ is
discontinuous, this metric fails to transition smoothly between
curvature $-h(\ell_1)$ and curvature $-1$ (in fact, $g_{t,0}(r)$ is
not even differentiable at $r=0$). To address this issue, we will find
values of $t$ near $h(\ell_1)$ and $\e$ near $0$ where the metric is
smooth and non-singular, and satisfies all the conditions of the
theorem.

\smallskip


First, note for any $\epsilon>0$ and any $0<t<1$, 
the functions $\fte$ and $\gte$ define a hyperbolic metric near $\bdy
V$.  On the interval $(-\epsilon/2, 0]$, $\kte(r)$ is identically $1$,
hence the differential equations satisfied by $\fte$ and $\gte$ are
solved by $\fte(r) = \ell_1 e^r$ and $\gte(r) = \ell_2 e^r$.  Thus in
the collar neighborhood $(-\epsilon/2, 0]$ of $0$, setting $f=\fte$
and $g=\gte$ in (\ref{eqn:riemann-metric}) gives the metric desired
near $\bdy V$.

The regularity of $\fte$ allows us to find a metric that is
non-singular. Recall that the cone angle along the core of $V$ will be
$2\pi$ whenever $f'(r_0)=2\pi$ (where $r_0$ is the root of $\fte(r)$).


\begin{lemma}\label{lemma:fprime=2pi} The roots of $\fte(r)$ have the
  following behavior:
\begin{itemize}
\item[(a)] For all $t \in (0,1)$ and $\e \geq 0$, $\fte(r)$ has a
  unique root $r_0(t,\e)$.
\item[(b)] The function $m(t,\e) := \fpte(r_0(t,\e))$ is continuous in
  $t$ and $\e$, and strictly decreasing in both variables.
\item[(c)] For every $t \in (0, h(\ell_1))$, there is a unique value
  $\e(t)>0$ such that $m(t, \e(t)) = 2\pi$.
\item[(d)] As $t \to h(\ell_1)$, $\e(t) \to 0$.
\end{itemize}
\end{lemma}

\begin{proof}
Parts (a) and (b) are proved in Lemma \ref{lemma:root-existence}. To
prove part (c), we study the explicit solution to the equation for
$\fte(r)$ when $\e = 0$. For all $r<0$, $f_{t,0}$ is given by the
simple differential equation
$$f_{t,0}''(r) = t\,f_{t,0}(r),$$
and the initial conditions $f_{t,0}(0) = \ell_1 = f_{t,0}'(0)$.  This
has solution:
\begin{eqnarray}\label{eq:f0-explicit}
  f_{t,0}(r) &=& \ell_1 \cosh \left( r\sqrt{t} \right) +
  \frac{\ell_1}{\sqrt{t}}\,\sinh \left( r\sqrt{t} \right) \nonumber \\
  &=& \frac{\ell_1\sqrt{1-t}}{\sqrt{t}}\sinh \left( \sqrt{t} \, 
    \left( r-r_0(t,0)\right) \right)
\end{eqnarray}
where $r_0(t,0) = -\tanh^{-1}(\sqrt{t})/\sqrt{t}$. Thus for all
$t\in(0, h(\ell_1)]$,
\begin{eqnarray*}
m(t,0) & = & f_{t,0}'(r_0(t,0)) \\
& = & \ell_1 \sqrt{1-t} \\
& \geq & \ell_1 \sqrt{ ( 2\pi / \ell_1 )^2 } \\
& = & 2\pi,
\end{eqnarray*}
with equality if and only if $t= h(\ell_1)$.

On the other hand, for all 
$\e > 2\ln(\ell_1/2\pi)$, $\: \fpte(-\e/2) = \ell_1 e^{-\e/2} < 2\pi$.
By its defining equation, $\fte(r)$ is concave up in $r$ when
$\fte(r)$ is positive, and concave down when negative.  Thus $m(t,\e)$
is the absolute minimum of $\fpte(r)$ over $\RR$.  Therefore
$$m(t,\e) < 2\pi \quad \mbox{whenever} \quad \e > 2\ln(\ell_1/2\pi).$$
By the intermediate value theorem, we can conclude that for all $t \in
(0, h(\ell_1))$, there is a value $\e(t)>0$ such that $m(t, \e(t)) =
2\pi$. 
Furthermore, by part (b), $m(t,\e)$ is strictly decreasing in $\e$,
and therefore $\e(t)$ is unique.

By part (b), we know that $m(t,\e)$ depends continuously on $t$ and
$\e$. Thus $\e(t)$ depends continuously on $t$. As a result, as $t \to
h(\ell_1)$, $\e(t) \to \e(h(\ell_1))$. Since we have already computed
that $m(h(\ell_1), 0) = 2\pi$, it follows that $\e(h(\ell_1)) = 0$,
completing the proof.
\end{proof}

From now on, we require that $t \in (0, h(\ell_1))$, and restrict our
attention to the functions $f_t := f_{t,\e(t)}$ and $g_t := g_{t,\e(t)}$
that give a non-singular Riemannian metric $\tau(t)$ on the solid
torus $V$. It remains to check the curvature and volume estimates for
this metric.

\begin{lemma}\label{lemma:curv-estimate}
  Fix a value of $t$ such that $\zeta h(\ell_1) \leq t < h(\ell_1)$.
  Then the Riemannian metric $\tau(t)$ defined by $f_t$ and $g_t$ has
  all sectional curvatures bounded above by $-\zeta h(\ell_1)$.
\end{lemma}

\begin{proof}
  We will actually prove the sectional curvatures of $\tau(t)$ are
  bounded above by $-t$.  Bleiler and Hodgson computed that these
  sectional curvatures are convex combinations of 
  $$\frac{-f''}{f}, \quad \frac{-g''}{g}, \quad \mbox{and} \quad
\frac{-f'\cdot g'}{f \cdot g}.$$
  By equations (\ref{eqn:def-f}) and
  (\ref{eqn:def-g}), we have
$$-\frac{f''_t(r)}{f_t(r)} = -k_{t,\e(t)}(r) \in [-1, -t] ,
 \quad -\frac{f'_t(r) \, g'_t(r)}{f_t(r) \, g_t(r)} =
-k_{t,\e(t)}(r) \in [-1, -t].$$

As for $g''_t(r)/g_t(r)$, we differentiate both sides of equation
(\ref{eqn:def-g}) to obtain
$$\frac{g''_t}{g_t} -
\left(\frac{g'_t}{g_t}\right)^2 \: = \: k_{t,\e(t)} -
\left(\frac{f_t}{f'_t}\right)^2
\frac{f''_t}{f_t}\,k_{t,\e(t)} +
\frac{f_t}{f'_t}\, k'_{t,\e(t)}\, ,$$
which simplifies, using equations (\ref{eqn:def-f}) and
(\ref{eqn:def-g}), to
$$\frac{g''_t}{g_t} \: = \: k_{t,\e(t)} +
\frac{f_t}{f'_t}\, k'_{t,\e(t)} \, .$$

Since $1 \geq k_{t,\e(t)} \geq t$ and all other terms are nonnegative
(because $f_t$ and $k_{t,\e(t)}$ are both increasing), $-g''_t/g_t
\leq -t$.
\end{proof}

\begin{lemma}\label{lemma:tube-volume}
  Let $t$ vary in the interval $(\zeta h(\ell_1), \, h(\ell_1))$, and
  define the \linebreak Riemannian metric $\tau(t)$ by the functions
  $f_t$ and $g_t$. Then
$$\lim_{t \to h(\ell_1)} \vol(V, \tau(t)) \: = \: \frac{\ell_1\ell_2}{2} \: = \: \half \, \area(\bdy
V).$$ 
\end{lemma}

\begin{proof}
By equation (\ref{eqn:riemann-metric}),
$$\vol(V, \tau(t)) = \int_{r_0(t, \e(t))}^0 f_{t,\e(t)}(r) \,
g_{t,\e(t)}(r) \, dr.$$ 
Let $\tlim := h(\ell_1)$. By Lemma \ref{lemma:fprime=2pi}, as $t \to
\tlim$, $\e(t) \to 0$. Furthermore, by Theorems \ref{thm:f-regular}
and \ref{thm:g-regular}, the functions $f_{t,\e(t)}$ and $g_{t,\e(t)}$
converge uniformly to $f_{\tlim, 0}$ and $g_{\tlim, 0}$, respectively.
Theorem \ref{thm:f-regular} also implies that
$$r_0(t,\e) := f^{-1}_{t,\e}(0)$$
is continuous in $t$ and $\e$. Thus, as $t \to \tlim$, $\, r_0(t,\e)
\to r_0(\tlim,0)$.
 
By equation (\ref{eq:f0-explicit}), we know that for $r<0$,
$$f_{\tlim, 0}(r) \: = \:  \frac{\ell_1\sqrt{1-\tlim}} {\sqrt{\tlim}}
\, \sinh \left( \sqrt{\tlim} \, ( r-r_0 ) \right),$$ 
where  $r_0 = -\tanh^{-1}(\sqrt{\tlim})/\sqrt{\tlim}$. 

Similarly, when $t=\tlim$ and $\e=0$, the differential equation for
$\gte$ has solution
$$g_{\tlim, 0}(r) \: = \:  \ell_2 \sqrt{1-\tlim} \, \cosh \left(
  \sqrt{\tlim} \, ( r-r_0 ) \right).$$ 
Thus we may compute:
$$\lim_ {t \to \tlim} \! \vol(V, \tau(t)) 
\; = \;  \lim_ {t \to \tlim}  
  \int_{r_0(t, \e(t))}^0 f_{t,\e(t)}(r) \, g_{t,\e(t)}(r) \, dr \hspace{1.5in} \textcolor{white}{.}$$
  \vspace{-3ex}
  {\setlength\arraycolsep{3pt}
\begin{eqnarray*}
\quad &=& \int_{r_0(\tlim, 0)}^0 f_{\tlim, 0}(r) \, g_{\tlim, 0}(r) \, dr \\
&=& \int_{r_0}^0 \!\ell_1 \ell_2 \frac{(1\!-\!\tlim)}{\sqrt{\tlim}}  
  \sinh \! \left( \! \sqrt{\tlim} \, ( r\!-\!r_0 ) \right)
  \cosh \!\left( \! \sqrt{\tlim} \, ( r\!-\!r_0 ) \right)  dr \\
&=& \left[ \ell_1 \ell_2 \, \frac{(1-\tlim)}{2\, \tlim} \, 
  \sinh^2 \! \left(\sqrt{\tlim} \, ( r-r_0 ) \right) \right]_{r_0}^0 \\ 
&=& \frac{\ell_1 \ell_2}{2} \cdot \frac{1-\tlim}{\tlim} \cdot 
  \sinh^2 \! \left( \tanh^{-1}(\sqrt{\tlim}) \right) \\
&=& \frac{\ell_1 \ell_2}{2} \cdot \frac{1-\tlim}{\tlim} \cdot 
  \frac{\tlim}{1-\tlim} \\
&=& \frac{\ell_1 \ell_2}{2}  \, .
\end{eqnarray*}
}
This completes the proof of Lemma \ref{lemma:tube-volume}.
\end{proof}

We are now ready to complete the proof of Theorem
\ref{thm:neg-curved-solid-torus}.  By Lemma \ref{lemma:fprime=2pi}, if
we select any $t < h(\ell_1)$ and $\e = \e(t) > 0$, we get a
non-singular metric satisfying conclusion (a) of the theorem.
By Lemma \ref{lemma:curv-estimate}, conclusion (b) is satisfied if we
ensure that $t$ is between $\zeta h(\ell_1)$ and $h(\ell_1)$. Finally,
by Lemma \ref{lemma:tube-volume}, if we select $t$ near enough to
$h(\ell_1)$, we will have $\vol(V) \geq \tfrac{\zeta}{2} \area(\bdy
V)$,
satisfying conclusion (c).
\end{proof}

\subsection{Negatively curved metrics on a 3--manifold}\label{sec:vol-comparison}

By applying Theorem \ref{thm:neg-curved-solid-torus} to several cusps
of a cusped manifold $M$, we obtain a negatively curved metric on a
Dehn filling of $M$.

\begin{theorem}\label{thm:neg-curved-manifold}
  Let $M$ be a complete, finite--volume hyperbolic manifold with
  cusps.  Suppose $C_1, \dots, C_k$ are disjoint horoball
  neighborhoods of some (possibly all) of the cusps.  Let $s_1, \dots,
  s_k$ be slopes on $\partial C_1, \dots, \partial C_k$, each with
  length greater than $2\pi$.  Denote the minimal slope length by
  $\lmin$.
  Let $S$ be the set of all Riemannian metrics on $M(s_1, \dots, s_k)$
  whose sectional curvatures lie in an interval $[-a,-1]$ for some
  constant $a \geq 1$. Then $S$ is non-empty, and
  $$\sup_{\sigma \in S} \, \vol(M(s_1, \dots, s_k), \sigma) \: \geq \:
  \left( h(\lmin) \right)^{3/2} \vol(M).$$
\end{theorem}

\begin{proof}
Fix an arbitrary constant $\zeta \in (0,1)$. We will replace each cusp
$C_i$ by a solid torus $V_i$ whose meridian is $s_i$.  Theorem
\ref{thm:neg-curved-solid-torus} guarantees the existence of a smooth
Riemannian metric $\tau_i$ on $V_i$, satisfying the following
properties:
\begin{itemize}
\item The sectional curvatures on $V_i$ are all at most 
$$
     - \zeta \, h(\ell(s_i)) \leq - \zeta \, h(\lmin).$$
\item $\displaystyle{\vol(V_i, \tau_i) \geq \tfrac{1}{2} \zeta \,
    \area(\bdy C_i) = \zeta \, \vol(C_i)}$.
\end{itemize}

Furthermore, in a neighborhood of each torus $\bdy C_i$, the metric
$\tau_i$ agrees with the hyperbolic metric on $M$.  Thus we may cut
out the cusps $C_1, \ldots, C_k$ and glue in the solid tori $V_1,
\ldots, V_k$, obtaining a smooth Riemannian metric $\tau$ on the
filled manifold $M(s_1, \dots, s_k)$, satisfying the following
properties:
\begin{itemize}
\item The sectional curvatures of $\tau$ are bounded above by $- \zeta
  \, h(\lmin)$, and below by some constant.  The lower bound comes
  from the fact that the solid tori $V_1, \dots, V_k$ are compact, and
  $\tau$ has constant curvature $-1$ on $M\setminus \cup_{i=1}^k V_i$.
\item $\vol(M(s_1, \dots, s_k), \tau) \geq \vol(M \setminus
\cup_{i=1}^k C_i) + \zeta\sum_{i=1}^k \vol(C_i)$ \newline
$\textcolor{white}{.} \hspace{1.3in}  \geq  \zeta \, \vol(M)$.
\end{itemize}

Now, we
would like our metric to have sectional curvatures bounded above by
$-1$.  Note the definition of sectional curvature implies that if we
rescale the metric $\tau$ to be $x\,\tau$, then all sectional
curvatures are multiplied by $x^{-2}$.  Thus we rescale $\tau$ to be
$\sigma = \sqrt{\zeta\,h(\lmin)}\,\tau$. 
This, in turn, rescales the volume by a factor $x^3 =
(\sqrt{\zeta\,h(\lmin)})^3$.  
Thus under the rescaled metric:
\begin{itemize}
\item The sectional curvatures of $\sigma$ lie in $[-a, -1]$ for
 some $a\geq 1$.
\item $\vol(M(s_1, \dots, s_k), \sigma) \geq \zeta^{5/2} \,
  (h(\lmin))^{3/2} \, \vol(M)$.
\end{itemize}

Thus we have found a metric $\sigma$ that lies in the set $S$. Now,
because $\zeta \in (0,1)$ was arbitrary, we can conclude that
  $$\sup_{\sigma \in S} \,  \vol(M(s_1, \dots, s_k), \sigma) \: \geq
  \: \left( h(\lmin) \right)^{3/2} \vol(M).$$ 
  
  \vspace{-4.5ex}
\end{proof}

To complete the proof of Theorem \ref{thm:vol-change}, suppose that
the manifold $N=M(s_1, \dots, s_k)$ admits a complete hyperbolic
metric $\sigma_{\rm hyp}$. (Since we have already proved that $N$
admits a negatively curved metric $\sigma$, the geometrization
conjecture implies that $N$ will indeed be hyperbolic.) Now, we
compare the volumes of these metrics via the the following theorem of
Boland, Connell, and Souto \cite{boland-connell-souto}, stated here in
a special case.

\begin{theorem}[\cite{boland-connell-souto}] \label{thm:BCS}
  Let $\sigma$ and $\sigma'$ be two complete,
  finite--volume Riemannian metrics on the same $3$--manifold $N$.
  Suppose that all sectional curvatures of $\sigma$ lie in the
  interval $[-1,1]$ and all sectional curvatures of $\sigma'$ lie in
  the interval $[-a,-1]$ for some constant $a \geq 1$. Then
  $$\vol(N, \sigma) \geq \vol(N, \sigma'),$$
  with equality if and only if both metrics are hyperbolic.
\end{theorem}

\begin{remark}
  When $N$ is a closed manifold, this theorem was originally proved by
  Besson, Courtois, and Gallot \cite{bcg}. In fact, it is quite likely
  that their proof would apply in our setting, because the negatively
  curved metrics that we construct all have constant curvature on the
  remaining cusps of $N$.
\end{remark}

\begin{proof}[Proof of Theorem \ref{thm:vol-change}]

By Theorem \ref{thm:neg-curved-manifold}, we know that $N=M(s_1,
\dots, s_k)$ admits a non-empty set $S$ of Riemannian metrics whose
sectional curvatures lie in an interval $[-a,-1]$.  By Theorem
\ref{thm:BCS}, the hyperbolic metric $\sigma_{\rm hyp}$ uniquely
maximizes volume over the set $S$. Thus, by putting together the
statements of the two theorems, we get:
$$\vol(N, \sigma_{\rm hyp}) \: = \: \max_{\sigma \in S} \,  \vol(N,
\sigma) \: \geq \: \left( h(\lmin) \right)^{3/2} \vol(M).$$
  
  \vspace{-4ex}
\end{proof}

\subsection{How sharp is Theorem \ref{thm:vol-change}?}\label{sec:numerics}
We will attempt to answer this question in two ways. For long slopes,
we compare the volume estimate of Theorem \ref{thm:vol-change} to the
asymptotic formula proved by Neumann and Zagier \cite{neumann-zagier}.
For medium--length slopes, we present the results of numerical
experiments conducted using SnapPea.

To compare asymptotic estimates, we restrict our attention to the case
when $M$ has exactly one cusp. Let $C$ be a maximal horoball
neighborhood of the cusp, let $s$ be a slope on $\bdy C$, and let
$$\Delta V \: := \: \vol(M) - \vol(M(s)).$$

With this notation, Neumann and Zagier \cite{neumann-zagier} proved
that as \linebreak $\ell(s) \to \infty$,
$$\Delta V \: \approx \: \frac{\pi^2 \, \area(\bdy C)}{\ell(s)^2} 
\: = \: \frac{2 \pi^2 \, \vol(C)}{\ell(s)^2} \, .$$ Meanwhile, by
expanding the Taylor series for $(1-x)^{3/2}$, we see that Theorem
\ref{thm:vol-change} implies
$$\Delta V \: \leq \: \frac{3}{2} \left( \frac{2\pi}{\ell(s)}
\right)^2 \vol(M) \: = \: \frac{6 \pi^2 \, \vol(M)}{\ell(s)^2} \, .$$

Thus, as $\ell(s) \to \infty$, Theorem \ref{thm:vol-change}
overestimates the change in volume by a factor of $3\vol(M)/\vol(C)$.
The quantity $\vol(C)/\vol(M)$ is known as the \emph{cusp density} of
$M$. B\"or\"oczky \cite{boroczky} has proved that the cusp density of
a hyperbolic manifold is at most $0.8533$. There is no known lower
bound on the cusp density; out of the approximately 5,000 orientable
cusped manifolds in the SnapPea census, exactly six have density less
than $0.45$. These numbers suggest that for most small manifolds,
Theorem \ref{thm:vol-change} overestimates the asymptotic change in
volume by a constant factor between $3.5$ and $7$.

\begin{figure}
\includegraphics[width=2.25in]{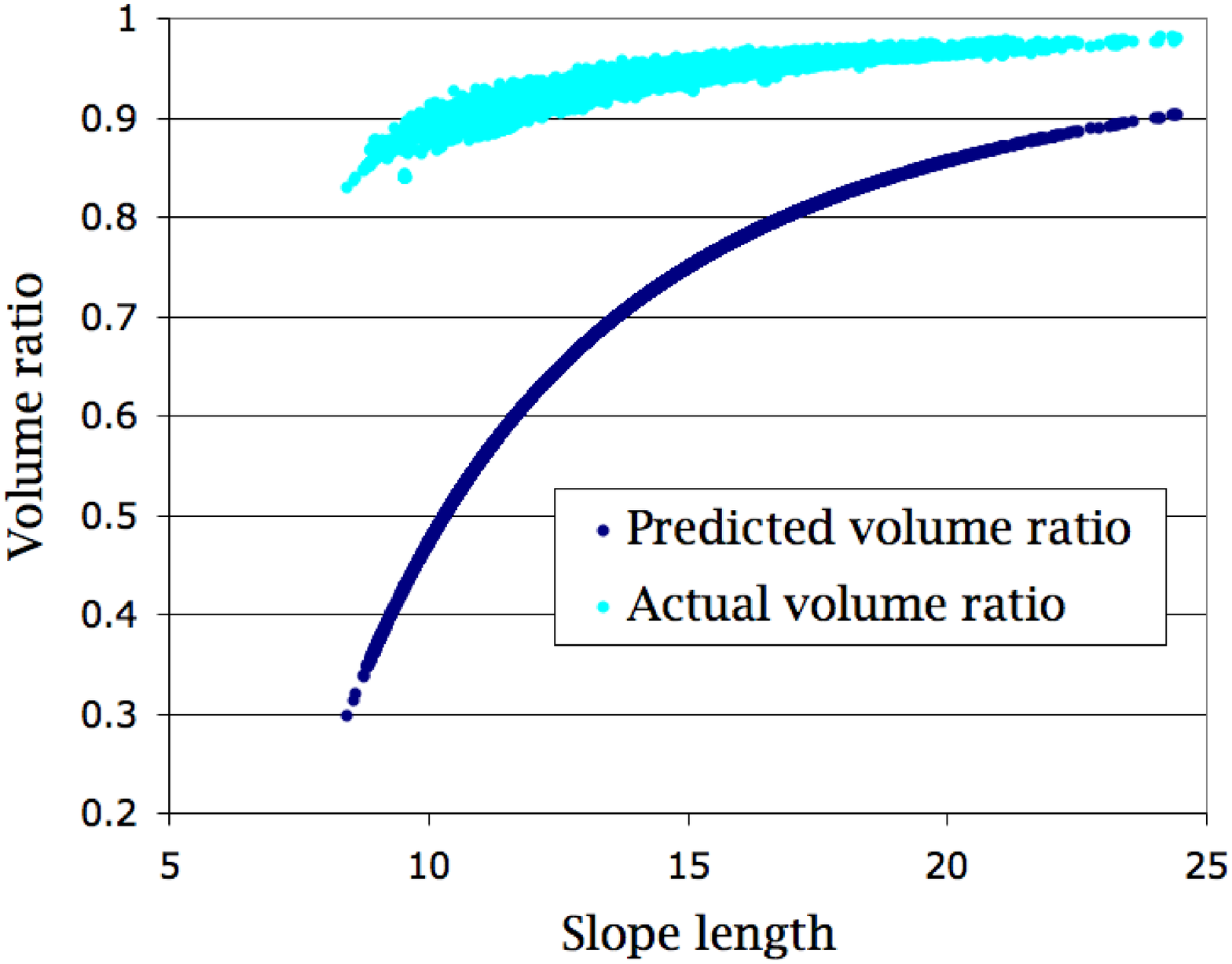}
\hspace{0.08in}
\includegraphics[width=2.25in]{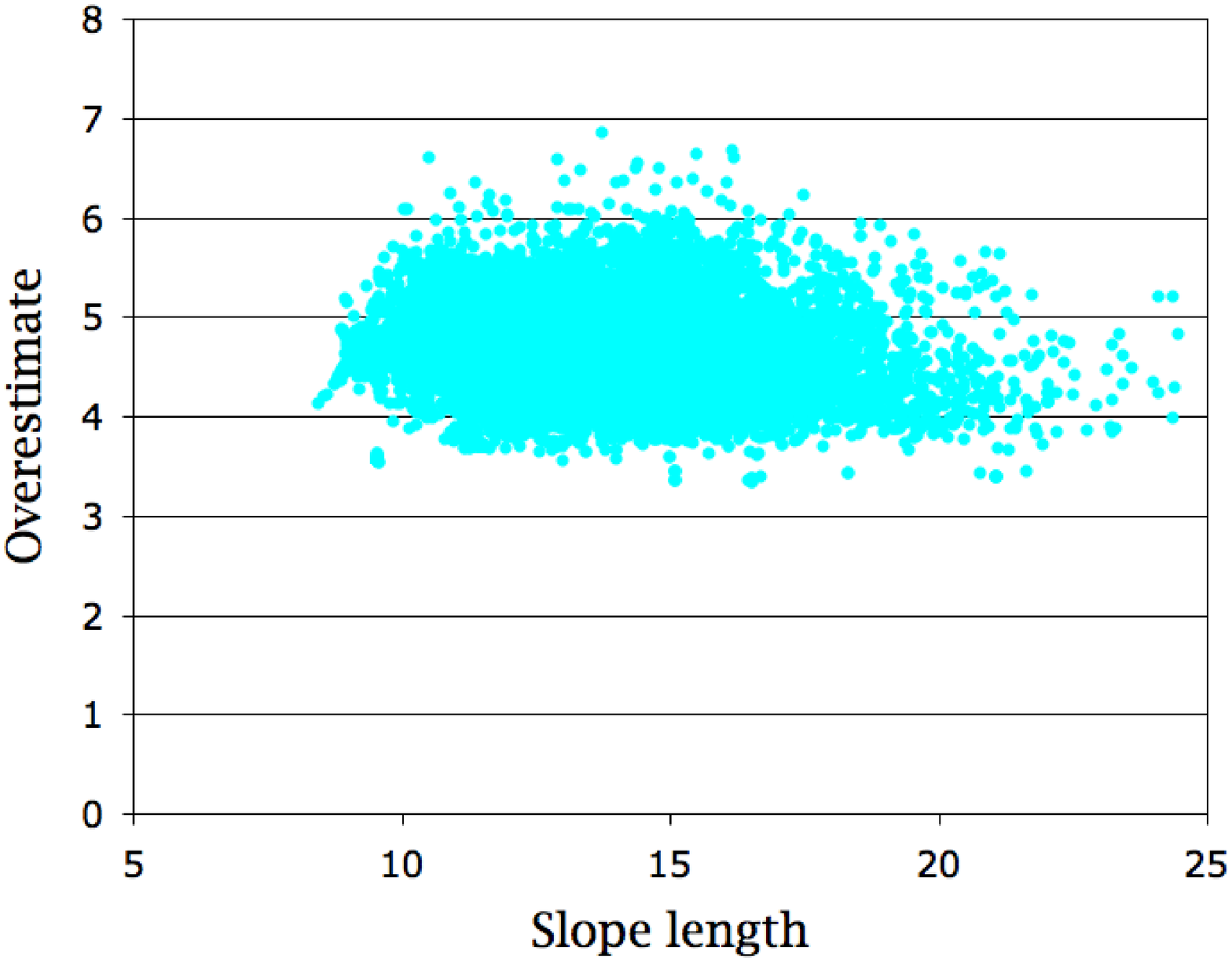}
\caption{The change in volume for medium--length slopes.}
\label{fig:numerical-data}
\end{figure}

For medium--length slopes, we also tested the estimate of Theorem
\ref{thm:vol-change} on over 14,000 manifold--slope pairs from the
SnapPea census. The results are plotted in Figure
\ref{fig:numerical-data}. In the graph on the left, the dark curve
represents the estimate of Theorem \ref{thm:vol-change}, while the
lighter point cloud represents the actual ratio $\vol(M(s))/\vol(M)$.
In the graph on the right, one can see that for all the manifolds and
slopes tested, Theorem \ref{thm:vol-change} overestimates the change
in volume by a factor between $3$ and $7$.

\section{Volumes of knots, links, and their fillings \label{sec:knots}}
In this section, we apply Theorem \ref{thm:vol-change} to hyperbolic link
complements in $S^3$, their Dehn fillings, and branched covers of $S^3$
over hyperbolic links.

\subsection{ Volumes of link complements}
To prove Theorem \ref{thm:link-volume}, we express a link $K$ as a
Dehn filling of another link $L$.

\begin{figure}[ht]
\psfrag{K}{$K$}
\psfrag{J}{$J$}
\psfrag{L}{$L$}
\psfrag{p}{$L'$}
\begin{center}

\includegraphics{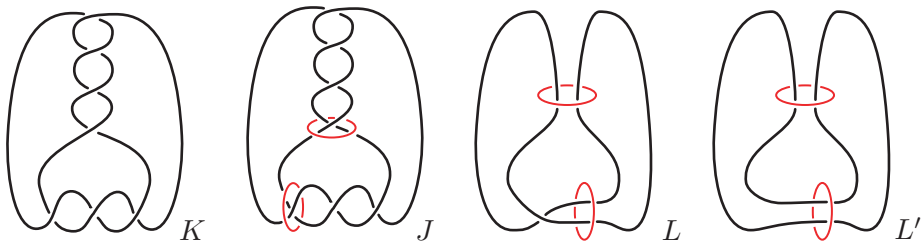}
\caption{An augmented link $L$ is constructed by adding a \emph{crossing
    circle} around each twist region of $D(K)$, then removing full twists.}
\label{fig:augment}
\end{center}
\end{figure}

Let $D(K)$ be a prime,
twist--reduced diagram of a link $K$ (see Section
\ref{sec:volume-results} for definitions). For every twist region of
$D(K)$, we add an extra link component, called a \emph{crossing
  circle}, that wraps around the two strands of the twist region. The result is a new link $J$. (See
Figure \ref{fig:augment}.)  Now, the manifold $S^3 \setminus J$ is homeomorphic to $S^3 \setminus L$, 
where $L$ is obtained by removing all full twists (pairs of crossings) from the twist regions of $J$.
This \emph{augmented link}
$L$ has the property that $K$ can be recovered by Dehn filling the
crossing circles of $L$. Similarly, every Dehn filling of $K$ can be
expressed as a filling of $L$.

The advantage of this construction is that the augmented link $L$ has
a simple geometry that allows for very explicit estimates.  

An estimate on volumes given by an estimate of cusp volume was given in
\cite{purcell:volume}.  Here we are able to improve that estimate.

\begin{prop}

\label{prop:aug-volume}
  Let $D(K)$ be a prime, twist--reduced diagram with at least two
  twist regions. Then the corresponding augmented link $L$ is
  hyperbolic, and
$$\vol(S^3 \setminus L) \; \geq \; 2\,v_8\,(\twist -1),$$
where $v_8 = 3.66386...$ is the volume of a regular ideal octahedron.
If $K$ is a two-bridge link, this inequality is an equality.
\end{prop}

\begin{proof}
The hyperbolicity of $S^3 \setminus L$ is a consequence of work of Adams
\cite{adams:aug}.  See also Purcell \cite{purcell:cusps}.  

To estimate the volume of $S^3 \setminus L$, we simplify the link $L$ even further, by removing all remaining single crossings from the twist regions of $L$. The resulting \emph{flat augmented link} $L'$ has the same volume as $L$, by the work of Adams \cite{adams:3-punct}. (See Figure \ref{fig:augment}.)  This link $L'$ is preserved by a reflection in the projection plane.
Thus the projection plane is isotopic to a totally
geodesic surface in $S^3 \setminus L'$.

Cut the manifold $S^3 \setminus L'$ along the projection plane.  The result is 
two hyperbolic manifolds $M$ and $M'$ with totally geodesic
boundary.  Since $M$ and $M'$ are interchanged by the reflection of
$S^3 \setminus L'$, they have the same volume.  Moreover, the volume of $S^3 \setminus L$
is given by the sum of the volumes of $M$ and $M'$.

Note that the manifold $M$ is a ball with a tube drilled out for each
crossing circle.  Hence it is topologically a handlebody with genus $\twist$.  Miyamoto showed that if $N$ is a hyperbolic $3$--manifold with totally geodesic
boundary, then $\vol(N) \geq -v_8\chi(N)$. (See \cite[Proposition 1.1 and Lemma 4.1]{miyamoto}.)  We apply
this result to $M$, and find 
$$\vol(S^3 \setminus L) \; = \; 2\, \vol(M) \; \geq \;  -2\, v_8 \, \chi(M) \; = \;  2\, v_8\,(\twist-1).$$  

Finally, when $D(K)$ is a standard diagram of a two-bridge link, it is well-known that the augmented link $L$ is obtained by gluing together $(\twist -1)$ copies of the Borromean rings, each of which has volume $2v_8$. (See, for example, Futer and Gu\'eritaud \cite[Theorem B.3]{gf-twobridge}.) Thus, for two-bridge links, $\vol(S^3 \setminus L) = 2v_8\,(\twist -1)$, making our estimate sharp.
\end{proof}

In fact, Proposition \ref{prop:aug-volume} is sharp for many additional large classes
of knot and link diagrams.

To recover $K$ from $L$, one must perform Dehn filling along the
crossing circles. Thus we need to estimate the
lengths of those slopes. To obtain information about Dehn fillings
of $K$, we also estimate the lengths of non-trivial (that is,
non-meridional) slopes on the components of $L$ that come from strands
of $K$.

\begin{prop}[Theorem 3.10 of \cite{futer-purcell}]\label{prop:aug-slopes}
  Let $K = \cup_{j=1}^{m} K_j$ be a link in $S^3$ with a prime,
  twist--reduced diagram $D(K)$. Suppose that $D(K)$ contains twist
  regions $R_1, \ldots, R_n$ ($n \geq 2$) and that twist region $R_i$
  contains $a_i$ crossings. For each component $K_j$, let $n_j$ be the
  number of twist regions visited by $K_j$, counted with multiplicity;
  and let $s_j$ be a non-trivial Dehn filling slope.

  Then the Dehn filling on $K$ with these slopes corresponds to a
  filling on the augmented link $L$. Furthermore, then there exists a
  choice of disjoint cusp neighborhoods in $S^3 \setminus L$, such
  that the slopes have the following lengths:

\begin{enumerate}
\item For a component $K_j$ of $K$, the slope $s_j$ has length at
  least $n_j$.

\item For a crossing circle $C_i$, the slope has length at least
  $\sqrt{a_i^2 + 1}$.
\end{enumerate}
\end{prop}

We now have enough information to prove Theorem \ref{thm:link-volume}.

\medskip

\noindent{\bf Theorem \ref{thm:link-volume}.}
{\it  Let $K \subset S^3$ be a link with a prime, twist--reduced diagram
  $D(K)$. Assume that $D(K)$ has $\twist \geq 2$ twist regions, and
  that each region contains at least $7$ crossings. Then $K$ is a
  hyperbolic link satisfying
  $$0.70735  \; (\twist-1) \; < \; \vol(S^3 \setminus K) \; < \; 10\,
  v_3 \, (\twist - 1), $$ where $v_3 \approx 1.0149$ is the volume of
  a regular ideal tetrahedron.  }

\begin{proof}
  The conclusion that $K$ is hyperbolic was proved
  by Futer and Purcell \cite[Theorem 1.4]{futer-purcell}, relying on
  W.\ Thurston's hyperbolization theorem \cite{thurston-survey}. (In
  fact, $6$ crossings per twist region suffice to show $K$ is
  hyperbolic.)  The upper bound on volume is due to Agol and D.\
  Thurston \cite{lackenby:alt-volume}.

  To prove the lower bound, we apply Theorem \ref{thm:vol-change} to
  the augmented link $L$. Since every twist region has at least $7$
  crossings, by Proposition \ref{prop:aug-slopes} the slope on each
  crossing circle will be at least $\sqrt{7^2+1} = 5\sqrt{2} >2\pi$.

  Thus, by Theorem \ref{thm:vol-change},
\begin{eqnarray*}
\vol(S^3\setminus K) &\geq&
\left(1-\left(\frac{2\pi}{5\sqrt{2}}\right)^2\right)^{3/2}
2 \: v_8 \: (\twist-1)\\
& = & 0.70735... \; (\twist-1).
\end{eqnarray*}

\vspace{-4ex}
\end{proof}

\begin{remark}
In the proof of Theorem \ref{thm:link-volume}, we used the
fact that to insert $7$ crossings into a
twist region, one fills along a slope of length at least $\sqrt{7^2+1} = 5\sqrt{2}$.  In fact,
if we require $8$ crossings per twist region, we may replace
$5\sqrt{2}$ with $\sqrt{8^2+1} = \sqrt{65}$, and the lower bound
improves to $1.8028\,(\twist-1)$. As the number of required crossings increases, the estimate becomes better still. In the case of $8$ crossings, our estimate is similar to the lower
bound for alternating links due to Lackenby \cite{lackenby:alt-volume}
and Agol, Storm, and Thurston \cite{ast-guts}, which is
$1.83(\twist-2)$.  Their estimate is known to be sharp for the
Borromean rings.
\end{remark}

\subsection{Dehn fillings and branched covers}

Under a slightly stronger diagrammatic condition than that of Theorem
\ref{thm:link-volume}, we can show that the combinatorics of a link
$K$ determines the volumes of all of its non-trivial fillings, up to
an explicit and bounded constant.

\begin{theorem}\label{thm:knot-surg-volume}
Let $K$ be a link in $S^3$ with a prime, twist--reduced diagram
  $D(K)$. Suppose that every twist region of $D(K)$ contains at least
  $7$ crossings and each component of $K$ passes through at least $7$
  twist regions (counted with multiplicity).  Let $N$ be a manifold
  obtained by a non-trivial Dehn filling of some (possibly all)
  components of $K$, which satisfies geometrization. Then $N$ is
  hyperbolic, and
  $$0.62768 \;  (\twist-1) \; < \; \vol(N) \; < \; 10\, v_3 \, (\twist - 1). $$
\end{theorem}

Note that if $K$ is a knot, a diagram with $4$ or more twist regions
and $7$ or more crossings per region satisfies the hypotheses of
Theorem \ref{thm:knot-surg-volume}.

The conclusion that every non-trivial filling of $K$ is hyperbolic was
first proved by Futer and Purcell \cite[Theorem 1.7]{futer-purcell},
modulo the geometrization conjecture. In fact, 6 crossings per twist
region suffice.

\begin{proof}
  To prove that $N$ is hyperbolic and compute the lower bound on
  volume, we once again apply Theorem \ref{thm:vol-change} to the
  augmented link $L$.  We know that every non-trivial filling of $K$
  can be realized as a filling of $L$. By Proposition
  \ref{prop:aug-slopes}, every slope on a strand of $K$ will have
  length at least $7$, and every slope on a crossing circle will have
  length at least $5\sqrt{2}$. Thus, by Theorem \ref{thm:vol-change},
  $N$ is hyperbolic and
\begin{eqnarray*}
\vol(N) &\geq & \left(1-\left(\frac{2\pi}{7}\right)^2\right)^{3/2}
\;2\;v_8\; (\twist-1) \\
& = & 0.62768... \;(\twist-1)
\end{eqnarray*}

For the upper bound, note that volume goes down under Dehn filling
(see Thurston \cite{thurston-notes}). Thus, by Theorem
\ref{thm:link-volume}, $\vol(N) < 10\, v_3 \, (\twist - 1)$.
\end{proof}

Theorem \ref{thm:vol-change} also applies to Dehn fillings of
arbitrary hyperbolic knots.

\begin{theorem} \label{thm:pq-volume} Let $N$ be a hyperbolic manifold
  obtained by $p/q$--Dehn surgery along a hyperbolic knot $K$ in
  $S^3$, where $\abs{q} \geq 12$.  Then
$$ \vol(N) \: > \:
  \left(1-\frac{127}{q^2}\right)^{3/2} \vol(S^3\setminus K).$$
\end{theorem}

\begin{proof}
  Let $C$ be a maximal cusp of $S^3 \setminus K$. Let $m$ and $s$ be
  Euclidean geodesics on $\bdy C$ that represent the meridian of $K$
  and the slope $p/q$, respectively. Let $\theta$ be the angle between
  these geodesics. Then
\begin{equation}\label{eq:cusp-area}
  \abs{q} \cdot \area(\bdy C) \: = \:  \ell(m) \, \ell(s) \,
  \sin(\theta) \: \leq \: \ell(m) \, \ell(s). 
\end{equation}

We can use equation (\ref{eq:cusp-area}) to estimate $\ell(s)$. By a
theorem of Cao and Meyerhoff \cite[Theorem 5.9]{cao-meyerhoff},
$\area(\bdy C) \geq 3.35$. Furthermore, by the 6--Theorem of Agol and
Lackenby \cite{agol:6theorem, lackenby:surgery}, surgery along a slope
of length more than $6$ yields a manifold with infinite fundamental
group, which cannot be $S^3$. Thus $\ell(m) \leq 6$. Combining these
results with equation (\ref{eq:cusp-area}) gives
\begin{equation}\label{eq:pq-estimate}
 \ell(s) \: \geq \:  \abs{q} \cdot 3.35 / 6.
\end{equation}
In particular, when $\abs{q} \geq 12$, $\ell(s) > 2\pi$. Plugging
inequality (\ref{eq:pq-estimate}) into Theorem \ref{thm:vol-change}
gives
$$ \! \! \! \vol(N) \: \geq \:
\left(\! 1- \left( \frac{6 \cdot 2\pi}{3.35 q} \right)^{\!\!2\,}
\right)^{\!3/2} \! \vol(S^3\setminus K) \: > \:
\left(\!1-\frac{127}{q^2}\right)^{\!3/2} \! \vol(S^3\setminus K).$$

\vspace{-3ex}
\end{proof}

We conclude the section with an application to branched covers.
Recall that the cyclic $p$--fold cover of a hyperbolic knot complement
$S^3\setminus K$ is a hyperbolic 3--manifold $X_p$ with torus
boundary. The meridian $m$ of $K$ lifts to a slope $m_p$ on $\bdy
X_p$. Then the $p$--fold branched cover of $S^3$ over $K$, denoted
$M_p$, is obtained by Dehn filling $\partial X_p$ along the slope
$m_p$.

\begin{theorem}
If $p \geq 4$, the branched cover $M_p$ is hyperbolic. For all $p
\geq 7$, we have
\begin{equation}\label{eq:branched-all}
\left(1-\frac{4 \pi^{2}}{p^2}\right)^{3/2}
\vol(S^3\setminus K) \: \leq \: {\vol(M_p)\over p}\: < \:
\vol(S^3\setminus K).
\end{equation}
If $K$ is not the figure--8 or $5_2$ knot and $p \geq 6$, the estimate improves to
\begin{equation}\label{eq:branched-better}
\left(1-\frac{2 \sqrt{2} \, \pi^{2}}{p^2}\right)^{3/2}
\vol(S^3\setminus K) \: \leq \: {\vol(M_p)\over p}\: < \:
\vol(S^3\setminus K).
\end{equation}
\end{theorem}
 
\begin{proof} 
The fact that $M_p$ is hyperbolic for $p\geq 4$ is a well--known consequence of the
orbifold theorem  (see e.g. \cite[Corollary 1.26]{chk}).

As for the volume estimate, the hyperbolic metric of $S^3 \setminus K$
lifts to a hyperbolic metric on $X_p$, implying that $\vol(X_p)=p \,
\vol(S^3\setminus K)$.  Furthermore, because a maximal cusp of $S^3
\setminus K$ lifts to a maximal cusp of $X_p$, we have
$\ell(m_p)=p\,\ell(m)$. To estimate the volume of $M_p$, we need to
estimate the length of $m$.

Adams has shown that every hyperbolic knot in $S^3$ has meridian of
length at least $1$ \cite{adams:waist}. Thus, in the cyclic cover
$X_p$, $\ell(m_p) \geq p$. In particular, when $p \geq 7$, we have
$\ell(m_p) > 2\pi$. Plugging $\ell(m_p) \geq p$ into Theorem
\ref{thm:vol-change} 
proves the estimate of equation (\ref{eq:branched-all}).

Adams has also proved that apart from the figure--8 and $5_2$
knots, every hyperbolic knot in $S^3$ has meridian of length at least
$2^{1/4}$ \cite{adams:waist2}.  Thus, in the cyclic cover
$X_p$, we have $\ell(m_p) \geq 2^{1/4} p$, proving equation (\ref{eq:branched-better}).
\end{proof}

\begin{remark}
  Numerical experiments with SnapPea confirm that the $p$--fold
  branched covers over the figure--8 and $5_2$ knots also satisfy
  equation (\ref{eq:branched-better}) when $6 \leq p \leq 1000$. The
  complements of these knots admit ideal triangulations consisting
  (respectively) of two and three tetrahedra, with simple gluing
  equations. Thus one can probably employ the methods of Neumann and
  Zagier \cite{neumann-zagier} to rigorously prove equation
  (\ref{eq:branched-better}) for these two knots. Most of the details
  of the figure--8 case are worked out in \cite[Section
  6]{neumann-zagier}.
\end{remark}

\section{Twist number and the {J}ones polynomial \label{sec:jones}}
In this section, we prove Theorem \ref{thm:twist-jones}. The proof has
three main steps. The first step, due to Stoimenow
\cite{stoimenow:coeff}, expresses the coefficients of the Jones
polynomial in terms of the combinatorics of the graphs $\GA$ and
$\GB$, defined in Definition \ref{def:ab-graphs}. The second and third
steps relate the combinatorics of the graphs to upper and lower bounds on the twist number of a
diagram.

\subsection {Reduced graphs and polynomial coefficients} 

\begin{define} \label{def:reduced-graphs} Let $D$ be an connected link
  diagram, with associated graphs $\GA$, $\GB$, as in Definition
  \ref{def:ab-graphs}.  The multiplicity of an edge $e$ of $\GA$ or
  $\GB$ is the number of edges that have their endpoints on the same
  pair of vertices as $e$. Let $\GRA$ denote the graph obtained from
  $\GA$ by removing multiple edges connected to the same
  pair of vertices. We will refer to $\GRA$ as the \emph{reduced
    $A$--graph} associated to $D$.  Similarly, the reduced $B$--graph
  $\GRB$ is obtained by removing multiple edges connected to the same
  pair of vertices.
  
  Let $v_A(D)$, $e'_A(D)$ (resp. $v_B(D)$, $e'_B(D)$) denote the
  number of vertices and edges of $\GRA$ (resp. $\GRB$).  When there
  is no danger of confusion we will omit $D$ from the notation above
  to write $v_A:=v_A(D)$, $v_B:=v_B(D)$, $e'_A:=e'_A(D)$ and
  $e'_B:=e'_B(D)$, and so on.
\end{define}

\begin{prop}[Stoimenow]\label{prop:oxs} 
  For a link diagram $D$, let
  $$\langle D \rangle = \alpha A^{m}+ \beta A^{m-4}+ \gamma A^{m-8}+
  \ldots + \gamma' A^{k+8}+ \beta' A^{k+4}+ \alpha' A^k$$
  denote the
  Kauffman bracket of $D$, so that $m$ (resp. $k$)
  is  the highest
  (resp. lowest) power in $A$. If $D$ is connected and $A$--adequate,
  then
$$\abs{\beta}= e'_A(D)- v_A(D)+1.$$
Similarly, if $D$ is connected and $B$--adequate, then
$$\abs{\beta'} = e'_B(D)- v_B(D)+1.$$
\end{prop}

Note it is well known that for an $A$--adequate diagram, $|\alpha|=1$,
and for a $B$--adequate diagram, $|\alpha'|=1$.

\begin{proof}
  For proofs of these statements, see the papers of Stoimenow
  \cite[Proposition 3.1]{stoimenow:coeff} or Dasbach and Lin
  \cite[Theorem 2.4]{dasbach-lin:head-tail}.
\end{proof}

To obtain the Jones polynomial $J_K(t)$ from the Kauffman bracket
$\langle D \rangle$, one multiplies $\langle D \rangle$ by a power of
$\minus A$ and sets $t:=A^4$. Thus the absolute values of the coefficients remain
the same. This gives the following immediate corollary:

\begin{corollary}\label{cor:coeff-expression}
  Let $D$ be an adequate diagram of a link $K$. Let $\beta$ and
  $\beta'$ be the second and next-to-last coefficients of $J_K(t)$.
  Then
$$\abs{\beta} + \abs{\beta'} \; = \; e'_A+e'_B- v_A-v_B + 2.$$
\end{corollary}

Given Corollary \ref{cor:coeff-expression}, we can complete the proof
of Theorem \ref{thm:twist-jones} by estimating the quantity
$e'_A(D)+e'_B(D)- v_A(D)-v_B(D) + 2$ in terms of $\twist$.

\subsection{Long and short resolutions}

\begin{figure}[ht]
\begin{center}
\input{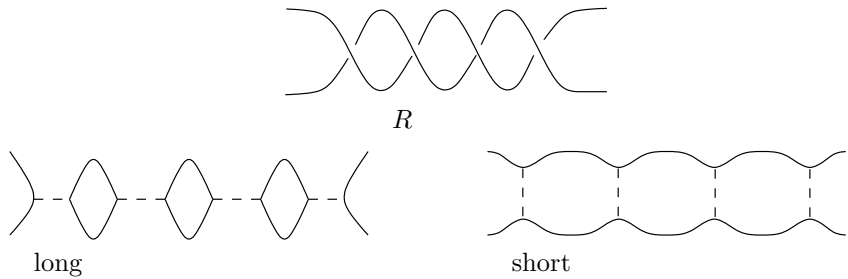}
\end{center}
\caption{Resolutions of a twist region $R$.}\label{fig:twist-resolutions}
\end{figure}

\begin{define}\label{def:long-short}
  Let $D$ be a diagram, and let $R$ be a twist region of $D$
  containing $c_R > 1$ crossings. One of the graphs associated to $D$,
  say $\GA$, will inherit $c_R -1$ vertices from the $c_R -1$ bigons
  contained in $R$. We say that this is the \emph{long resolution} of
  the twist region $R$. The other graph, say $\GB$, contains $c_R$
  parallel edges (only one of which survives in $\GRB$). This is the
  \emph{short resolution} of $R$. See Figure
  \ref{fig:twist-resolutions}.
  
  When a twist region $R$ contains a single crossing, there is no
  natural way to choose the short and long resolutions. For such a
  twist region, we say that both resolutions are short.
\end{define}

In order to count the vertices and edges of $\GRA$ and $\GRB$, we
regroup them into short and long resolutions.

\begin{define}\label{def:regrouping}
  Recall from Definition \ref{def:ab-graphs} that every vertex of
  $\GA$ and $\GB$ (and thus of $\GRA$ and $\GRB$) comes from a component
  of one of the $A$-- or $B$--resolutions of the diagram $D$. We say
  that a vertex adjacent to exactly $2$ edges of $\GA$ or $\GB$ is a
  \emph{bigon vertex}; these vertices correspond to bigons in twist
  regions of $D$. Let $\vbigon$ be the total number of bigon vertices
  in $\GA$ and $\GB$, and let $\vlarge$ be the total number of
  remaining, non-bigon vertices of $\GA$ and $\GB$.
  
  In a similar vein, let $\eshort$ (resp. $\elong$) be the total
  number of edges of $\GA$ and $\GB$ coming from short (resp. long)
  resolutions of twist regions. Observe that an edge comes from a long
  resolution if and only if it is adjacent to at least one bigon
  vertex.  Thus, when a pair of vertices is connected by multiple
  edges, if neither vertex is a bigon vertex, those edges all short.
  In any other case, those edges are all long. As a result, we can
  think of every edge of $\GRA$ and $\GRB$ as either short or long,
  and define $\ershort$ and $\erlong$ accordingly.
\end{define}

An immediate consequence of this definition is that
$$v_A+v_B = \vbigon + \vlarge \quad \mbox{and} \quad e'_A+e'_B =
\erlong + \ershort.$$

We are now ready to prove one direction of Theorem
\ref{thm:twist-jones}.

\begin{prop}\label{thm:reduced-at-most} 
  Let $D$ be an adequate diagram of a link $K$. Let $\beta$ and
  $\beta'$ be the second and next-to-last coefficients of the Jones
  polynomial $J_K$. Then
  $$\abs{\beta} + \abs{\beta'} \; = \; e'_A+e'_B- v_A- v_B + 2 \; \leq
  \; 2\, \twist.$$
\end{prop}

\begin{proof}
  Suppose that an adequate diagram $D$ has $c:=c(D)$ crossings and
  $t:=\twist$ twist regions.  Given Corollary
  \ref{cor:coeff-expression} and Definition \ref{def:regrouping}, it
  suffices to estimate the quantities $\vbigon$, $\vlarge$, $\erlong$,
  and $\ershort$ in terms of $c$ and $t$.
  
  In a twist region $R$ containing $c_R$ crossings, there are $c_R-1$
  bigons. Thus $\vbigon = c-t$. Notice that in both the $A$-- and
  $B$--resolutions of $D$, at least one circle passes through multiple
  twist regions. Thus each of $\GA$ and $\GB$ contains at least one
  non-bigon vertex, and $\vlarge \geq 2$. In every twist region, all
  the edges of the short resolution get identified to a single edge in
  either $\GRA$ or $\GRB$. Thus $\ershort \leq t$. Meanwhile, since
  each crossing has at most one long resolution, $\erlong \leq \elong
  \leq c$. Putting these facts together, we get

\vspace{.07in}
\begin{tabular}{rcll}
$\abs{\beta} + \abs{\beta'}$&$=$& $e'_A+e'_B- v_A- v_B + 2$, & (Corollary
  \ref{cor:coeff-expression})  \vspace{.05in}\\ 
&$=$&  $\ershort + \erlong - \vbigon - \vlarge +2,
$ & (Definition \ref{def:regrouping}) \vspace{.05in} \\ 
&$\leq$& $t \quad + \quad c \quad - \;\;  (c-t) \;\; -
  \; 2 \;\; + 2$ &  \vspace{.05in}\\ 
&$=$& $2t$.&
\end{tabular}
\vspace{.07in}

We note that the adequacy of $D$ is only needed to apply Corollary
\ref{cor:coeff-expression}. The remainder of the proof works for any
connected diagram $D$.
\end{proof}

\subsection{Estimates from Turaev surfaces} 

To obtain a lower bound on $\abs{\beta} + \abs{\beta'}$, we engage in
the detailed study of a \emph{Turaev surface} associated to the $A$--
and $B$--resolutions of a diagram. The construction of this surface
was first described by Cromwell \cite{cromwell-book}, building on work
of Turaev \cite{turaev}.

\begin{figure}[ht]
\psfrag{sa}{$s_A$}
\psfrag{sb}{$s_B$}
\psfrag{g}{$\Gamma$}
\begin{center}
\includegraphics{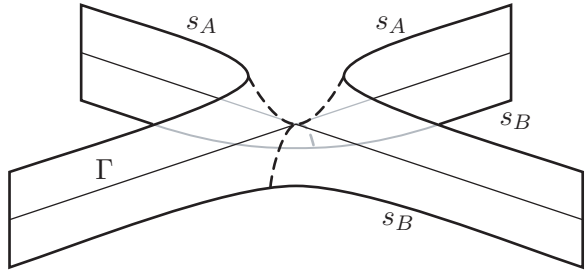}
\caption{Near each crossing of the diagram, a saddle surface interpolates 
  between circles of $s_A(D)$ and circles of $s_B(D)$. The edges of
  $\GA$ and $\GB$ can be seen as gradient lines at the saddle.}
\label{fig:saddle}
\end{center}
\end{figure}

Let $\Gamma \subset S^2$ be the planar, 4--valent graph of the link
diagram $D$.  Thicken the projection plane to a slab $S^2 \cross
[\minus 1, 1]$, so that $\Gamma$ lies in $S^2 \cross \{0\}$. Outside a
neighborhood of the vertices (crossings), our surface will intersect
this slab in $\Gamma \cross [\minus 1, 1]$. In the neighborhood of
each vertex, we insert a saddle, positioned so that the boundary
circles on $S^2 \cross \{1\}$ are the
components
of the $A$--resolution $s_A(D)$, and the boundary circles on $S^2
\cross \{\minus 1\}$ are the components of $s_B(D)$. (See Figure
\ref{fig:saddle}.) Then, we cap off each circle with a disk, obtaining
an unknotted closed surface $F(D)$.

In the special case when $D$ is an alternating diagram, each circle of
$s_A(D)$ or $s_B(D)$ follows the boundary of a region in the
projection plane. Thus, for alternating diagrams, the surface $F(D)$
is exactly the projection sphere $S^2$. For general diagrams, it is
still the case that the knot or link has an alternating projection to
$F(D)$ \cite[Lemma 4.4]{dessin-jones}.

Furthermore, the construction of $F(D)$ endows it with a natural
cellulation, whose $1$--skeleton is the graph $\Gamma$ and whose
$2$--cells correspond to circles of $s_A(D)$ or $s_B(D)$, hence to
vertices of $\GA$ or $\GB$. These $2$--cells admit a natural
checkerboard coloring, in which the regions corresponding to the
vertices of $\GA$ are white and the regions corresponding to $\GB$ are
shaded. The graph $\GA$ (resp. $\GB$) can be embedded in $F(D)$ as the
adjacency graph of white (resp. shaded) regions.

\begin{define}\label{def:geography}
  Let $D$ be a diagram in which every twist region has at least $2$
  crossings (hence, at least one bigon). Then we may modify the
  $4$--valent graph $\Gamma \subset F(D)$, by collapsing the chain of
  bigons in each twist region to a single \emph{red edge}. The result
  is a tri-valent graph $P \subset F(D)$, in which exactly one
  edge at each vertex is colored red. (See Figure
  \ref{fig:surface-graphs}.)

\begin{figure}[ht]
\psfrag{k}{$K$}
\psfrag{g}{$\Gamma$}
\psfrag{d}{$P$}
\psfrag{f}{$\Phi$}
\psfrag{ar}{$\Rightarrow$}
\begin{center}
  \includegraphics{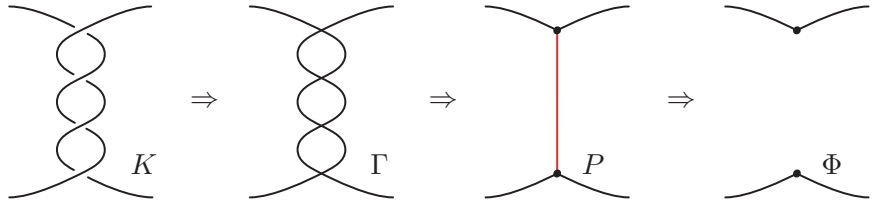}
\caption{The construction of the graphs $P$ and $\Phi$ from a
  knot diagram. The entire construction takes place on the surface
  $F(D)$.} 
\label{fig:surface-graphs}
\end{center}
\end{figure}

If we remove all the red edges of $P$, we obtain a di-valent
graph $\Phi$. In other words, $\Phi$ is a union of disjoint simple
closed curves. We call the closures of regions in the complement of $P$ the
\emph{provinces} of $F(D)$, and the closures of regions in the complement of
$\Phi$ the \emph{countries} of $F(D)$.
\end{define}

Every province of $F(D)$ comes from a non-bigon region of $F(D)
\setminus \Gamma$, and corresponds to
an n-gon
vertex of
$\GA$ or $\GB$. Thus each province is a disk. On the other hand, the
countries of $F(D)$ need not be simply connected.

The subdivision of $F(D)$ into countries allows us to partition the
twist regions of $D$ and the short edges of $\GA$ and $\GB$. Every
twist region of $D$ corresponds to a red edge that lies in some
country $N$ of $F(D)$. Similarly, every short edge of $\GA$ or $\GB$
connects two vertices that belong to the same country $N$. (Compare
Figure \ref{fig:twist-resolutions} with Figure
\ref{fig:surface-graphs}.) Thus we may define $\twistcountry$ to be
the number of twist regions belonging to $N$, and $\eshort(N)$ to be
the number of short edges belonging to $N$. In a similar fashion, we
may define $\ershort(N)$ by removing the duplicate edges of $\GA$ or
$\GB$ that belong to $N$.

\begin{lemma}\label{lemma:short-count}
  Let $N$ be a country of $F(D)$. Then
  $$\ershort(N) \geq \twistcountry + \chi(N) - 1.$$
\end{lemma}

\begin{proof}
The country $N$ is constructed by taking a number of contractible
provinces and gluing them along disjoint segments (red edges) on the
boundary.
We claim if we cut $N$ along a well--chosen set of $1 - \chi(N)$ red
edges, it becomes a disk. This can be seen by considering the dual
graph to the red edges. Note the country deformation retracts to this
dual graph. A maximal spanning tree is obtained by removing
$1-\chi(N)$ edges, which correspond to red edges in $N$.
  
After this operation, there remain $\twistcountry + \chi(N) -1$ red
edges along which we did not cut. Call these the \emph{remnant red
  edges}.
The remnant edges are in one-to-one correspondence with a subset of
elements of $\ershort(N)$, given by selecting a short edge of $\GA$ or
$\GB$ from the corresponding twist region.  So $\ershort(N) \geq
\twistcountry + \chi(N) -1$.
\end{proof}

To estimate $\ershort(D)$ more globally, we need a bound on the number
of countries.

\begin{lemma}\label{lemma:country-count}
  Let $D$ be an adequate diagram, in which every twist region contains
  at least $2$ crossings. Let $n(D)$ be the number of countries in the
  surface $F(D)$. Then
$$n(D) \leq \frac{2}{3} \, \twist + 1.$$
\end{lemma}

\begin{proof}
Recall, from Definition \ref{def:geography}, that every national
border is a component of the graph $\Phi$.
Let $\abs{\Phi}$ denote the number of components of the graph $\Phi$.
Thus $n(D) \leq \abs{\Phi} +1$. Observe as well that the graph $P$ has
exactly two vertices for every twist region of $D$ (these are the two
endpoints of the red edge constructed from the twist region). Thus we
compare the number of components of $\Phi$ to the number of vertices
of $P$.
  
Let $\varphi$ be one simple closed curve contained in $\Phi$. We will
count the number of vertices of $P$ that lie on $\varphi$.  There are
two straightforward cases:
  
\smallskip

\emph{\underline{Case 1: $\varphi$ bounds a one--province country.}}
This province cannot be a \linebreak monogon, since monogons cannot
occur in an adequate diagram. It also cannot be a bigon, because the
bigon would have been collapsed in the construction of the graph $P$.
Thus $\varphi$ contains at least $3$ vertices of $P$.
  
\smallskip

\emph{\underline{Case 2: $\varphi$ does not bound a one--province
    country.}} Then consider the \linebreak provinces that adjoin
$\varphi$.  The provinces of $F(D)$ are simply connected, so each side of
$\varphi$ must meet at least one provincial border (red edge).  In
fact, the hypothesis that $D$ is adequate implies that a province
cannot border on itself along a red edge (otherwise, an edge of $\GA$
or $\GB$ dual to this red edge would form a loop, violating Definition
\ref{def:ab-graphs}). Thus each side of $\varphi$ must meet at least
two provinces, so $\varphi$ must contain at least $4$ vertices of $P$.
  
\smallskip
  
In either case, each curve $\varphi \subset \Phi$
contains at least 3 vertices from twist regions.  Since each twist
region gives rise to two such vertices,
$\twist \geq \frac{3}{2}
\abs{\Phi}$. We can conclude that
  $$n(D) \; \leq \; \abs{\Phi} + 1 \; \leq \; \frac{2}{3} \, \twist
  +1.$$

\vspace{-4ex}
\end{proof}

We can now prove the remaining direction of Theorem
\ref{thm:twist-jones}.

\begin{theorem} \label{thm:reduced-at-least}
  Let $D$ be an adequate diagram of a link $K$, in which every twist
  region contains at least $3$ crossings. Let $\beta$ and $\beta'$ be
  the second and next-to-last coefficients of the Jones polynomial
  $J_K$. Then
  $$\abs{\beta} + \abs{\beta'} \; = \; e'_A+e'_B- v_A - v_B + 2 \;
  \geq \; \frac{\twist}{3} + 1.$$
\end{theorem}
 
\begin{proof} 
  Suppose that the diagram $D$ has $c:=c(D)$ crossings and $t:=\twist$
  twist regions.  Given Corollary \ref{cor:coeff-expression} and
  Definition \ref{def:regrouping}, it suffices to estimate the
  quantities $\vbigon$, $\vlarge$, $\erlong$, and $\ershort$ in terms
  of $c$ and $t$.
  
  In a twist region $R$ containing $c_R$ crossings, there are $c_R-1$
  bigons and $c_R$ long edges. Thus $\vbigon = c-t$ and $\elong = c$.
  When every twist region contains at least $3$ crossings, it is
  evident from Figure \ref{fig:twist-resolutions} that all long edges
  of $\GA$ and $\GB$ will survive in $\GRA$ and $\GRB$. Thus we can
  conclude that $\erlong = c$, giving us
\begin{equation}\label{eq:long}
\erlong - \vbigon = \twist.
\end{equation}

To estimate $\vlarge$ and $\ershort$, we compute the Euler
characteristic of $F(D)$. Recall that the tri-valent graph $P$
has $2t$ vertices (two for every red edge) and $3t$ edges (since every
third edge is red). The $2$--cells in the complement of $P$ are
provinces, one for every
n-gon
vertex of $\GA$ and $\GB$. Thus
\begin{equation}\label{eq:euler}
\chi(F(D))  \; = \; \vlarge -3\twist + 2\twist \; = \; \vlarge - \twist.
\end{equation}

Lemma \ref{lemma:short-count} tells us that $\ershort(N) \geq
\twistcountry + \chi(N) - 1$ for every country $N$ of $F(D)$. By
summing this over all countries, we get

\vspace{.07in}
\begin{tabular}{rcll}
$\ershort$ &$\geq$& $\twist + \chi(F(D)) - n(D)$ & \vspace{.05in}\\
&$=$& $\vlarge - n(D)$, & by Equation
  (\ref{eq:euler})\vspace{.05in}\\ 
&$\geq$& $\displaystyle{\vlarge - \frac{2}{3} \, \twist - 1}$ & by
  Lemma \ref{lemma:country-count}.\vspace{.05in}
\end{tabular}
\vspace{.07in}

\noindent Putting all of these results together gives

\vspace{.07in}
\begin{tabular}{rcll}
$\abs{\beta} + \abs{\beta'}$ &$=$& $e'_A+e'_B- v_A+v_B + 2$, &by Corollary
  \ref{cor:coeff-expression}\vspace{.05in}\\ 
&$=$& $(\erlong - \vbigon) + (\ershort - \vlarge) +2, \quad$ & by
  Definition \ref{def:regrouping}\vspace{.05in}\\ 
  &$\geq$& $\displaystyle{\twist \quad +\; \;  \left(- \frac{2}{3} \, \twist - 1
    \right) + 2},$ & \vspace{.05in} \\
  &$=$& $\displaystyle{\frac{\twist}{3} \, + \, 1}$.
\end{tabular}

\vspace{-3ex}
\end{proof}

\section{Families of differential equations \label{sec:diffeq}}

In the proof of Theorem \ref{thm:neg-curved-solid-torus}, we defined a
family of functions $\fte$ and $\gte$.  Our goal in this section is to
prove that $\fte$, $\fpte$, and $\gte$ depend continuously and
uniformly on $t$ and $\e$.

Let us recap the definitions.  For parameters $\e >0$ and $0 < t < 1$,
we began with a smooth bump function $\kte(r)$. 
This function has a precise definition, as follows:
$$\kte(r) \; := \; 
\left\{ 
\begin{array}{rl}
t & \quad \mbox{if} \quad r \leq - \epsilon, \\
\displaystyle{ t \, + \, (1-t) \, \frac{\int_0^{2+2r/\epsilon} z(u) \, du}
{\int_0^1 z(u) \, du} }
&\quad \mbox{if} \quad -\epsilon < r < -\epsilon/2,\\
1 & \quad \mbox{if} \quad r \geq - \epsilon/2,
\end{array}
\right.
$$
where $z(u) = \exp \left(-\frac{1}{u^2}-\frac{1}{(u-1)^2} \right)$.  (See Figure \ref{fig:bump-function} for a typical graph.)

By extension to $\e=0$, we defined $k_{t,0}(r)$ as a step function
whose value is $t$ for $r<0$ and $1$ for $r \geq 0$.

Given $\kte$, we defined $\fte$ and $\gte$ according to the
differential equations
\begin{equation}\label{eqn:def-fg}
\fppte(r) = \kte(r) \, \fte(r), \quad \quad
\frac{\gpte(r)}{\gte(r)} = \kte(r) \, \frac{\fte(r)}{\fpte(r)} \, ,
\end{equation}
with initial conditions $\fte(0)= \fpte(0) = \ell_1$ and $\gte(0) =
\ell_2$.

We will prove that $\fte$, $\fpte$, and $\gte$ depend continuously
and uniformly on $t$ and $\e$, even as $\e$ goes to $0$, when $\kte$
becomes discontinuous. Before we prove that statement, we need a
monotonicity result.

\begin{define}\label{def:r0}
  For any $t \in (0,1)$ and $\e \geq 0$, define
  $$r_0(t,\e) := \inf{s \in \RR : \fte(r) > 0 \mbox{ for all } r >
    s}.$$
  In other words, $r_0$ is either the largest root of $\fte$,
  or $-\infty$ if $\fte$ has no root. 
\end{define}

\begin{lemma}\label{lemma:r0}
  For any $t \in (0,1)$ and $\e \geq 0$, $\: \: r_0(t, \e) < \min \{ -\e/2, -1 \}$.
\end{lemma}

\begin{proof}
If $r_0(t,\e) = -\infty$, the result is trivially true. Thus we may assume that $r_0(t,\e)$ is a root of $\fte$.
Note that for $r \geq -\e/2$,  \linebreak $\fte(r) =\ell_1 e^r >  0$, and thus $r_0(t, \e) <-\e/2 \leq 0$.

To prove that $r_0(t,\e) < -1$, observe that equation (\ref{eqn:def-fg}) implies $\fte(r)$ is concave up on $(r_0(t,\e), 0]$. Thus, for all $r \in (r_0(t,\e), 0)$ we have
$$\fpte(r) \: < \: \fpte(0) \: = \: \ell_1.$$
Since the function $\fte$ must climb from height $0$ to height $\ell_1$ with slope less than $\ell_1$, it follows $r_0(t,\e) < -1$.
\end{proof}


\begin{lemma}\label{lemma:f-monotonic}
  The functions $\fte$ and $\fpte$ are monotonic in the parameters $t,
  \e$; and $\fte$ is also monotonic in $r$. More precisely:
\begin{itemize}
\item[(a)] If $0 \leq \e_1 < \e_2$ and $r \in [r_0(t,\e_2),\,0]$, then
  $f_{t,\e_1}(r) \leq f_{t,\e_2}(r)$ and $f'_{t,\e_1}(r) \geq
  f'_{t,\e_2}(r)$, with strict inequalities on $[r_0(t,\e_2), \,
  -\e_1/2)$.
\item[(b)] If $0 < t_1 <t_2 < 1$ and $r \in [r_0(t_2,\e),\,0]$, then
  $f_{t_1,\e}(r) \leq f_{t_2,\e}(r)$ and $f'_{t_1,\e}(r) \geq
  f'_{t_2,\e}(r)$, with strict inequalities on $[r_0(t_2,\e),\,
  -\e/2)$.
\item[(c)] For any $r \in \RR$, $\e \geq 0$, and $t \in (0,1)$,
  $\fpte(r) >0$.
\end{itemize}
\end{lemma}

\begin{proof}
The key observation for this proof is that the bump function $\kte(r)$
is increasing in both $t$ and $\e$. See Figure
\ref{fig:bump-function}.

For part (a), suppose that $0 \leq \e_1 < \e_2$. To compare $f$ and
$f'$ for these two values of $\e$, define a function $\phi(r) :=
f_{t,\e_2}(r) - f_{t,\e_1}(r)$. Then
\begin{eqnarray*}
\phi''(r) &=& f''_{t,\e_2}(r) - f''_{t,\e_1}(r) \\
&=& k_{t,\e_2}(r) \, f_{t,\e_2}(r) - k_{t,\e_1}(r) \, f_{t,\e_1}(r) \\ 
&\geq& k_{t,\e_1}(r) \, f_{t,\e_2}(r) - k_{t,\e_1}(r) \, f_{t,\e_1}(r)
  \quad \mbox{when } r \in [r_0(t,\e_2),\,0],\\ 
& & \quad \mbox{with a strict inequality for } r \in [r_0(t,\e_2),\,0]
  \cap (-\e_2, -\e_1/2)   \\ 
&=& k_{t,\e_1}(r) \, \phi(r).
\end{eqnarray*}
By definition, $0 \leq k_{t,\e_1}(r) \leq 1$.
Thus we obtain a differential inequality with certain nice properties.
By a result from analysis, whose proof we include as Theorem
\ref{thm:diff-inequal} in the Appendix, $\phi(r) \geq 0$ and $\phi'(r)
\leq 0$ for all $r \in [r_0(t,\e_2),\,0]$, with strict inequalities on
$[r_0(t,\e_2), \, -\e_1/2)$. Note that this interval is non-empty,
because by Lemma \ref{lemma:r0},
$$r_0(t,\e_2) \: < \:  -\e_2/2 \: < \: -\e_1/2.$$
This proves (a).

The proof of part (b) is very similar to (a), except this time we
define $\phi(r) := f_{t_2,\e}(r) - f_{t_1,\e}(r)$. An analogous
calculation then goes through.

For part (c), fix values of $\e \geq 0$ and $t \in (0,1)$. We want to
prove that $\fpte(r)>0$ for all $r$. If $r \geq -\e/2$, we have
already seen that $\fpte(r)=\ell_1 e^r >0$.

If $r_0(t,\e) \leq r < -\e/2$, we rely on part (a). That is, set $\e_1
= \e$ and $\e_2 = - 2r$.  Then
provided we can show $r \in [r_0(t, \minus 2r), 0]$,
part (a) implies
$$\fpte(r) \: \geq\: f'_{t, \minus 2r}(r) \: = \: \ell_1 e^r \: > \:
0.$$
If $r_0(t, \minus 2r)=-\infty$, then certainly $r$ is in the correct
range.  Otherwise, we know $r \geq r_0(t,\e)$ by assumption.  For any
$s \in [r_0(t, \minus 2r), 0]$, part (a) implies $\fte(s) \leq f_{t,
  \minus 2r}(s)$, so in particular, $\fte(r_0(t, \minus 2r)) \leq
f_{t, \minus 2r}(r_0(t, \minus 2r))=0$.  Thus $r_0(t,\e) \geq r_0(t,
\minus 2r)$.  So $r\geq r_0(t, \minus 2r)$ as desired.

Finally, for $r < r_0(t,\e)$, we note that when $\fte(r)<0$, the
function must also be concave down. Thus, since $\fpte(r) >0$ at the
root $r=r_0(t,\e)$, $\fpte(r)$ can only become more positive as $r$
moves further to the left.
\end{proof}

\begin{theorem}\label{thm:f-regular}
  Fix constants $\tlim \in (0,1)$ and $\elim \geq 0$. Then, as $(t, \e) \to
  (\tlim, \elim)$,
  $$\fte(r) \to f_{\tlim,\elim}(r) \quad \mbox{and} \quad \fpte(r) \to
  f_{\tlim,\elim}'(r),$$
  uniformly on compact sets. In particular, the
  functions $f$ and $f'$ are continuous in the three variables $(t,\e,
  r)$.
\end{theorem}

\begin{proof}
Observe that when $\e >0$, the function $\kte(r)$ is continuous in all
three variables $(t,\e, r)$. Thus, when $\elim >0$, the conclusion of
the theorem is a standard result in ODE theory (see, for example,
\cite{hsieh-sibuya}). We will therefore restrict our attention to the
case when $\elim =0$.

Fix an integer $n$ such that $\tlim \in (\frac{1}{n}, \frac{n-1}{n})$.
Now, suppose that $(t, \e)$ varies in the compact domain
$[\frac{1}{n}, \frac{n-1}{n}] \times [0, 1]$, and that $r$ varies in
the compact interval $[-n, n]$. We begin the argument by showing that
the values of $\fte(r)$ are uniformly bounded on this domain. By Lemma
\ref{lemma:f-monotonic}(c), $\fte(r)$ is strictly increasing, and so
attains its maximum value at $r=n$. Thus $\fte(n) = \ell_1 e^n$ is a
uniform upper bound.

For a lower bound on $\fte(r)$, we take a closer look at equation
(\ref{eqn:def-fg}). When $r \leq -\e$, the equation has the explicit
solution
\begin{equation}\label{eq:f-explicit}
\fte(r) \; =\;  c_1(t, \e) e^{r \sqrt{t}} + c_2(t, \e) e^{-r \sqrt{t}},
\end{equation}
where
\begin{equation}\label{eq:c1c2}
c_1(t, \e) \: = \: \frac{\fte(-\e)}{2 e^{-\e\sqrt{t}}} +
\frac{\fpte(-\e)}{2 \sqrt{t}e^{-\e\sqrt{t}}} \, , \quad 
c_2(t, \e) \: = \: \frac{\fte(-\e)}{2 e^{\e\sqrt{t}}} -
\frac{\fpte(-\e)}{2 \sqrt{t}e^{\e\sqrt{t}}}\, . 
\end{equation}
Now, when $\e \leq 1$, Lemma \ref{lemma:r0} says that $r_0(t, \e) < -\e$. Thus $\fte$ and $\fpte$ are
positive and increasing on $[-\e, 0]$, and both are bounded above by
$\ell_1$. Thus both $\fte(-\e)$ and $\fpte(-\e)$ must be in  the interval $(0, \ell_1]$. 
In particular, this implies that
\begin{equation}\label{eq:c1c2-estimate}
c_1(t, \e) > 0, \quad c_2(t,\e) > \frac{-\ell_1}{2 \sqrt{1/n}} \, .
\end{equation}
By putting together (\ref{eq:f-explicit}) and
(\ref{eq:c1c2-estimate}), we see that on the interval $[-n, n]$,
$$\fte(r) \: \geq \: \fte(-n) \: > \: c_2(t,\e) \, e^{n \sqrt{t}} \:
\geq \: \frac{-\ell_1 \sqrt{n}}{2} \, e^n.$$

We can conclude that when $t \in [\frac{1}{n}, \frac{n-1}{n}]$, $\e
\in [0, 1]$, and $r \in [-n, n]$, the family of functions $\fte(r)$
is uniformly bounded. Because $\kte(r)$ is also uniformly bounded (by
0 and 1), it follows that $\fppte(r)$ is uniformly bounded. By
integration, it follows that $\fpte(r)$ is uniformly bounded and
equicontinuous. Integrating again, we see that $\fte(r)$ is
equicontinuous. Finally, for any $\delta \in (0, n)$, $\kte(r)$ has
uniformly bounded derivative on $[-n, -\delta]$, which implies that
$\fppte(r)$ is equicontinuous on that interval.

Fix a number $\delta \in (0, n)$, and let $(t_i, \e_i)$ be a sequence
that converges to $(\tlim, 0)$. Then the Arzela--Ascoli theorem
implies that there is a continuous function $f_{\rm lim}$ on
$[-n, n]$,
twice differentiable on $[-n, -\delta]$, such that
$$f_{t_i, \e_i}(r) \to f_{\rm lim}(r), \quad f'_{t_i, \e_i}(r) \to
f'_{\rm lim}(r), \quad f''_{t_i, \e_i}(r) \to f''_{\rm lim}(r),$$
uniformly on $[-n, -\delta]$. In fact, $f_{t_i, \e_i}$ and $f'_{t_i,
  \e_i}$ converge uniformly on $[-n, n]$.
Furthermore, for all $\e_i < \delta$, $k_{t_i, \e_i}(r) = t_i$ on
$[-n, -\delta]$, and thus $k_{t_i, \e_i}(r)$ converges uniformly to
$\tlim$. Thus $f_{\rm lim}$ satisfies the differential equation
$$f''_{\rm lim}(r) = k_{\tlim, 0}(r) f_{\rm lim}(r),$$
for all $r \in
[-n, -\delta]$. Since $\delta$ was arbitrary, this equation is
satisfied for all $r \in [-n, 0)$. Since $f_{\rm lim}(0) = f'_{\rm
  lim}(0) = \ell_1$, $f_{\rm lim}$ is a solution to equation
(\ref{eqn:def-fg}), for $t=\tlim$ and $\e=0$. Therefore, by the
uniqueness of solutions, we can conclude that $f_{\rm lim}(r) =
f_{\tlim, 0}(r)$, for all $r \in [-n,n]$.
\end{proof}


\begin{lemma}\label{lemma:root-existence} The roots of $\fte(r)$ have
  the following behavior:
\begin{itemize}
\item[(a)] For all $t \in (0,1)$ and $\e \geq 0$, $\fte(r)$ has a
  unique root, equal to $r_0(t,\e)$.
\item[(b)] The function $m(t,\e) := \fpte(r_0(t,\e))$ is continuous in
  $t$ and $\e$, and strictly decreasing in both variables.
\end{itemize}
\end{lemma}

\begin{proof}
By Lemma \ref{lemma:f-monotonic}(c), $\fte(r)$ is strictly increasing
on $\RR$. Thus if a root exists, it will be unique. To prove the
existence of a root, we study the explicit formula for $\fte(r)$ on
the interval $(-\infty,-\e]$, given in equation (\ref{eq:f-explicit}).
As $r \to -\infty$, this equation is dominated by the term $c_2(t,\e)
\, e^{-r \sqrt{t}}$.  In particular, since $\fte$ is increasing on
$\RR$, we must have $c_2(t,\e) \leq 0$. We will show that, in fact, $c_2(t,\e)
<0$.

Suppose, for a contradiction, that $c_2(t,\e_1) = 0$ for some value
$\e_1$. Then, by equation (\ref{eq:f-explicit}),
$$f_{t,\e_1}(r) = c_1(t, \e_1) \, e^{r \sqrt{t}} \quad \mbox{on }
(-\infty, -\e_1],$$
for a positive constant $c_1(t,\e_1)$. Now, choose
a larger value $\e_2$. As $\e_2 \to \infty$, we have larger and larger
subsets of $(-\infty,0]$ on which
$$f_{t,\e_2}(r) = \ell_1 e^r.$$
Because $t<1$ and $e^r$ decays faster
than $e^{r \sqrt{t}}$ as $r\to -\infty$, there will be an $\e_2 \gg
\e_1$ and an $r \ll 0$ such that
$$f_{t,\e_2}(r) \: = \: \ell_1 e^r \: < \: c_1(t, \e_1) \, e^{r
  \sqrt{t}} \: = \: f_{t,\e_1}(r),$$
contradicting Lemma
\ref{lemma:f-monotonic}(a). Thus $c_2(t,\e) <0$ for all $t,\e$.

As a result, $\fte(r)$ approaches $-\infty$ as $r \to -\infty$, and
therefore has a root. By Definition \ref{def:r0}, this unique root is
equal to $r_0(t, \e)$.

To prove part (b), we once again use the fact that $\fte(r)$ is
continuous and strictly increasing. Thus it has a continuous inverse
$f^{-1}_{t,\e}$, such that $f^{-1}_{t,\e}(0) = r_0(t,\e)$. This allows
us to write
$$m(t,\e) = \fpte \circ f^{-1}_{t,\e}(0).$$
By Theorem
\ref{thm:f-regular}, both $\fpte$ and $f^{-1}_{t,\e}$ are continuous
in $t$ and $\e$; therefore, $m(t,\e)$ is continuous as well.

Now, fix starting values $\e_1, t_1$ of $\e$ and $t$. Then, for any
$\e_2 > \e_1$, Lemma \ref{lemma:f-monotonic}(a) implies that
$$
m(t_1,\e_1)= f'_{t_1,\e_1}(r_0(t_1,\e_1)) >
f'_{t_1,\e_2}(r_0(t_1,\e_1)).$$
Since $m(t_1,\e_2)$ is the absolute minimum of $f'_{t_1,\e_2}(r)$
over all of $\RR$ (because $\fte$ is concave up whenever $\fte$ is
positive, concave down when negative due to its defining equation
(\ref{eqn:def-fg})), we have
$$f'_{t_1,\e_2}(r_0(t_1,\e_1)) \geq m(t_1,\e_2).$$
So $m(t_1,\e_1) \geq m(t_1,\e_2)$.

Similarly, by Lemma \ref{lemma:f-monotonic}(b), $m(t_1,\e_1) > m(t_2,
\e_1)$ for $t_2>t_1$. Thus $m(t,\e)$ is strictly decreasing in both
$t$ and $\e$.
\end{proof}


We now turn our attention to the function $\gte(r)$. Its defining
equation (\ref{eqn:def-fg}) can be written as
\begin{equation}
\frac{d}{dr} \left( \ln \gte(r) \right) \: = \: 
\kte(r) \, \frac{\fte(r)}{\fpte(r)} \, ,
\label{eqn:redef-g}
\end{equation}
with initial condition $\gte(0) = \ell_2$.  Note that by Lemma
\ref{lemma:f-monotonic}(c), \linebreak $\fpte(r)>0$ for all $r$, so the
right-hand side is always well-defined.

\begin{theorem}\label{thm:g-regular}
  Fix constants $\tlim \in (0,1)$ and $\elim \geq 0$. Then, as $(t, \e) \to
  (\tlim, \elim)$,
$$\gte(r) \to g_{\tlim, \elim}(r),$$
uniformly on compact sets.
\end{theorem}

\begin{proof}
This proof follows the same outline as the proof of Theorem
\ref{thm:f-regular}. As in that proof, we restrict our attention to
the case when $\elim =0$, because the conclusion of the theorem is a
standard result for $\elim > 0$.

Fix an integer $n$ such that $\tlim \in (\frac{1}{n}, \frac{n-1}{n})$.
Now, suppose that $(t, \e)$ varies in the compact domain
$[\frac{1}{n}, \frac{n-1}{n}] \times [0, 1]$, and that $r$ varies in
the compact interval $[-n, n]$. We begin the argument by showing that
the right-hand side of equation (\ref{eqn:redef-g}) is uniformly
bounded on this domain. In the proof of Theorem \ref{thm:f-regular},
we have already shown that on this domain,
\begin{equation}\label{eq:bound1}
 \abs{ \kte(r) \, \fte(r) } \: \leq \:  \abs{\fte(r)} \: \leq \:
 \ell_1 \sqrt{n} \, e^n. 
\end{equation}
Also, because $m(t,\e)$ is the absolute minimum value of $\fpte(r)$
over all of $\RR$, and by Lemma \ref{lemma:root-existence}(b),
\begin{equation}\label{eq:bound2}
\fpte(r) \: \geq \: m(t, \e) \: \geq \: m(\tfrac{n-1}{n}, \, 1) \: >\: 0.
\end{equation}

Putting inequalities (\ref{eq:bound1}) and (\ref{eq:bound2}) together,
we get
\begin{equation}\label{eq:bound3}
\left| \frac{d}{dr} \left( \ln \gte(r) \right) \right| \: \leq \:
\frac{\ell_1 \sqrt{n} \, e^n}{m(\frac{n-1}{n}, \, 1)} \, .
\end{equation}

By integrating (\ref{eq:bound3}), we conclude that the family of
functions $\ln \gte(r)$ is uniformly bounded and equicontinuous. Also,
for any $\delta \in (0, n)$, $\frac{d}{dr} \ln \gte(r) $ is
equicontinuous on $[-n, -\delta]$. This follows by differentiating the
right-hand side of (\ref{eqn:redef-g}), because $\fte(r)$, $\fpte(r)$,
$\fppte(r)$, $\kte(r)$, and $k'_{t,\e}(r)$ are all uniformly bounded
on that interval, with $\fpte(r)$ bounded away from 0.

Fix a number $\delta \in (0, n)$, and let $(t_i, \e_i)$ be a sequence
that converges to $(\tlim, 0)$. Then the Arzela--Ascoli theorem
implies that there is a continuous function $g_{\rm lim}$,
differentiable on $[-n, -\delta]$, such that
$$\ln g_{t_i, \e_i}(r) \to \ln g_{\rm lim}(r), \quad \frac{d}{dr} \ln
g_{t_i, \e_i}(r) \to \frac{d}{dr} \ln g_{\rm lim}(r), $$
uniformly
on $[-n, -\delta]$. In fact, $\ln g_{t_i, \e_i}(r)$ converges
uniformly on $[-n, n]$; since this is a compact set, $g_{t_i,
  \e_i}(r)$ also converges uniformly to $g_{\rm lim}(r)$.

By letting $\delta$ approach $0$, we see that the function $g_{\rm
  lim}(r)$ satisfies the differential equation (\ref{eqn:redef-g}) for
$t = \tlim$ and $\e = 0$. Thus, by the uniqueness of solutions,
$g_{\rm lim}(r) = g_{\tlim,0}(r)$, as desired.
\end{proof}

\appendix
\section{Differential inequalities}


The following elementary result from real analysis is probably
well-known. However, since we could not find a reference, we include a
proof here.

\begin{theorem}\label{thm:diff-inequal}
  Let $I \subset \RR$ be
a closed
  interval that includes $0$. Let $\phi: I
  \to \RR$ be a $C^1$ function, such that $\phi''(x)$ exists for all
  $x \neq 0$. Suppose that $\phi$ satisfies the differential
  inequality
  $$\phi''(x) \geq k(x) \: \phi(x) \quad \mbox{for all } x \neq 0,$$
  where $0 \leq k(x) \leq 1$. Assume as well that $\phi(0) \geq 0$ and
  $\phi'(0)=0$.  Then
\begin{itemize}
\item[(a)] $\phi'(x) \geq 0$ for $x \geq 0$ and $\phi'(x) \leq 0$ for
  $x\leq 0$,
\item[(b)] $\phi(x) \geq 0$ and $\phi''(x) \geq 0$ for all $x$.
\end{itemize}
Furthermore,
\begin{itemize}
\item[(c)] If $\phi''(x_0)>0$ and $x_0 <0$, then $\phi(x)>0$ and
  $\phi'(x)<0$ for $x < x_0$.
\item[(d)] If $\phi''(x_0)>0$ and $x_0 >0$, then $\phi(x)>0$ and
  $\phi'(x)>0$ for $x > x_0$.
\end{itemize}
\end{theorem}

\noindent A key step of the proof is the following, slightly weaker statement.

\begin{lemma}\label{lemma:inductive}
  Let $I \subset \RR$ be a closed
interval that includes $0$. Let $\psi: I
  \to \RR$ be a $C^1$ function, such that $\psi''(x)$ exists for all
  $x \neq 0$. Suppose that $\psi$ satisfies the differential
  inequality
  $$\psi''(x) \geq k(x) \: \psi(x) \quad \mbox{for all } x \neq 0,$$
  where $0 \leq k(x) \leq 1$. Assume as well that $\psi(0) \geq 0$ and
  $\psi'(0) \geq 0$. Then
\begin{itemize}
\item[(a)] $\psi(x) \geq 0$, $\psi'(x) \geq 0$, and $\psi''(x) \geq 0$
  on $[0,1] \cap I$.
\item[(b)] If $\psi''(x_0) >0$ or $\psi'(x_0)>0$ for some $x_0 \in
  [0,1)$, then $\psi(x)>0$ and $\psi'(x)>0$ on $(x_0, 1] \cap I$.
\end{itemize}
\end{lemma}

\begin{proof}
  To prove (a), let $m= \min \{ \psi(x): x \in [0,1] \cap I \}$.
  Assume, for a contradiction, that $m<0$. Then, because $k(x) \leq 1$
  for all $x$, we have $\psi''(x) \geq m$, for all $x \in [0,1] \cap I
  $. Now, Taylor's theorem allows us to write
  $$\psi(x) = \psi(0) + \psi'(0) \, x + \half \, \psi''(x_0) \, x^2,
  \quad \mbox{for some } x_0 \in (0,x).$$
  We can estimate each of
  these terms. We have $\psi(0) \geq 0$ by hypothesis, $\psi'(0) \, x
  \geq 0$ because both parts of the product are non-negative, and
  $\psi''(x_0) \, x^2 \geq m$ because $x^2 \leq 1$. (Recall that we
  have assumed $m<0$.) Putting all of this together gives
  $$\psi(x) \geq 0 + 0 + \half \, m > m \quad \mbox{for all } x \in [0,1]
  \cap I,$$
  contradicting the assumption that $m$ was the minimum.
  
  As a result of this contradiction, $\psi(x) \geq 0$ on $[0,1] \cap
  I$. Thus, since $k(x) \geq 0$, we have $\psi''(x) \geq 0$ as well.
  Integration gives $\psi'(x) \geq 0$, completing the proof of (a).
  
  To prove (b), suppose first that $\psi''(x_0) >0$ for some $x_0 \in
  [0,1)$. Then $\psi'(x)$ is strictly increasing in a neighborhood of
  $x_0$.  Since we have already shown that $\psi''(x) \geq 0$ for all
  $x \in [0,1] \cap I$, we have $$\psi'(x) > \psi'(x_0) \geq 0, \quad
  \mbox{for all } x \in (x_0, 1] \cap I.$$
  By integration, we also have
$\psi(x) > 0$ on $(x_0, 1] \cap I$.
  
  Now, suppose that $\psi'(x_0) >0$ for some $x_0 \in [0,1)$. Since
  $\psi''(x) \geq 0$ for all $x$, this implies that $\psi'(x) > 0$ on
  $[x_0, 1] \cap I$. Then, integration gives $\psi(x) > 0$ on $(x_0,
  1] \cap I$, completing the proof.
\end{proof}

\begin{proof}[Proof of Theorem \ref{thm:diff-inequal}]
  We will apply Lemma \ref{lemma:inductive} inductively, many times.
  We assume without loss of generality that $I = \RR$; if $I
  \subsetneq \RR$, the only alteration required is to stop the
  inductive process once we get to the boundary of $I$.

  Applying the lemma to $\psi(x) := \phi(x)$ gives the conclusion of
  the theorem on the interval $[0,1]$.  That is, $\phi(x)$ satisfies
  the differential inequality of the lemma, and $\phi(0)\geq 0$, and
  $\phi'(0) =0$ so the lemma applies immediately.  This will imply
  that $\phi(1) \geq 0$ and $\phi'(1) \geq 0$, with strict
  inequalities if $\phi''(x_0) > 0$ for some $x_0 \in [0,1)$.
  
  Now, apply Lemma \ref{lemma:inductive} to the function $\psi(x) :=
  \phi(x+1)$.  
  The lemma applies because
  $\psi(0)=\phi(1)\geq 0$, and $\psi'(0)=\phi'(1)
  \geq 0$ (possibly with strict inequalities).  This gives the
  conclusion of the theorem on the interval $[1,2]$. Repeatedly
  applying the lemma in this way proves the theorem for all $x \geq
  0$.

  To prove the theorem for $x \leq 0$, we first apply Lemma
  \ref{lemma:inductive} to the function $\psi(x) := \phi(-x)$.  Note
  that $\psi(0)=\phi(0)\geq 0$, and $\psi'(0) = -\phi'(0) = 0$, so the
  lemma applies to this function.  Then we obtain the conclusion of
  the theorem on $[-1, 0]$. Note, in particular, that we now have
  $\phi(-1) \geq 0$ and $\phi'(-1) \leq 0$ (with strict inequalities
  if $\phi''(x_0) > 0$ for some $x_0 \in (-1,0]$).  Now, apply Lemma
  \ref{lemma:inductive} to $\psi(x) := \phi(-x-1)$, etc., to obtain the
  conclusion of the theorem for all $x \leq 0$.
\end{proof}

\bibliographystyle{hamsplain}

\bibliography{biblio.bib}
\end{document}